\documentclass[11pt]{amsart}
\usepackage{amsmath, amscd}
 
\theoremstyle{plain}
\newtheorem{thm}{Theorem}[section]
\newtheorem{lm}[thm]{Lemma}
\newtheorem{cor}[thm]{Corollary}
\newtheorem{pro}[thm]{Proposition}

\theoremstyle{definition}
\newtheorem{df}[thm]{Definition}

\def \R {\mathbb{R}}
\def \C {\mathbb{C}}
\def \Z {\mathbb{Z}}

\def \CB {\mathcal B}
\def \CE {\mathcal E}

\def \J {\mathcal J}
\def \CP {\mathcal P}

\def \CM {\mathcal M}
\def \hM {\hat{{\mathcal M}}}
\def \CR {\mathcal R}

\def \la {\langle}
\def \ra {\rangle}
\def \re {representation}
\def \a {\alpha}
\def \b {\beta}
\def \c {\chi}
\def \g {\gamma}

\def \d {\delta}

\def \lam {\lambda}

\def \n {\nabla}
\def \o {\omega}
\def \ov {\overline}
\def \p {\phi}
\def \r {\rho}
\def \s {\sigma}
\def \Si {\Sigma}
\def \t {\tau}
\def \th {\theta}
\def \ti {\tilde}

\def \O {\Omega}

\def \bd {\partial}
\def \ob {\overline{\partial}}
\def \pj {\overline{\partial}_J}

\def \ox {\otimes}
\def \op {\oplus}
\def \x {\times}
\def \ve {\varepsilon}

\begin{document}
 
\baselineskip.525cm
 
\title[Floer Homology, Composite Knots]
{The symplectic Floer homology of composite knots}
\thanks{Partially supported by NSF Grant DMS 9626166}
\author[Weiping Li]{Weiping Li}
\address{Department of Mathematics \newline
\hspace*{.175in} Oklahoma State University \newline
\hspace*{.175in}Stillwater, Oklahoma 74078-0613\newline
\hspace*{.175in}U. S. A.}
\email{wli@littlewood.math.okstate.edu}
\date{February 2, 1998}

\begin{abstract}
We develop a method of calculation for the symplectic Floer homology of
composite knots. The symplectic Floer homology of knots defined in
\cite{li} naturally admits an integer graded lifting, and it formulates a
filtration and induced spectral sequence. Such a spectral sequence
converges to the symplectic homology of knots in \cite{li}.
We show that there is another spectral
sequence which converges to the $\Z$-graded symplectic Floer homology 
for composite knots represented by braids.
\end{abstract}
\maketitle

\section{Introduction\label{Intro}}

For integral homology 3-spheres $Y$, Casson defined an integral invariant
which roughly counts the number of irreducible $SU(2)$ {\re}s of the
fundamental group $\pi_1(Y)$; Floer developed a $\Z_8$-graded instanton homology
theory based upon an application of Morse theory to the Chern-Simons functional
on the space ${\CB}_Y$ of equivalence classes of $SU(2)$ connections
on $Y$ (see \cite{fl1}). The Euler characteristic of the instanton Floer
homology is twice Casson's invariant. Lin \cite{lin} studied the special {\re}s
of knot group $\pi_1(S^3 \setminus K)$ into $SU(2)$ such that all meridians
of knot $K$ are represented by trace-zero matrices, and defined an invariant
$\lam_{CL}(K)$ for the knot $K$. 
The present author \cite{li} 
developed a symplectic Floer homology based upon the Atiyah 
conjecture that there is no difference between the instanton Floer homology
and the symplectic Floer homology for the integral homology 3-spheres.
Our symplectic Floer homology generalizes the invariant of Lin, and
its Euler characteristic is the negative of $\lam_{CL}(K)$.
For $SU(2)$ {\re}s of $\pi_1(S^3 \setminus K)$ with the trace of all
meridians fixed (not necessary zero), Cappell, Lee and Miller \cite{clm},
independently Herald \cite{he}, defined the equivariant knot signature
from the symplectic theory and the gauge theory points of view; the present
author extended the symplectic Floer homology in \cite{li} to the
general case which the Euler characteristic is the equivariant knot signature
in \cite{li2}.

In this paper, we study the Mayer-Vietoris principle for the symplectic
Floer homology defined in \cite{li, li2}. We only restrict to the trace-zero
case in \cite{li} (the general case in \cite{li2} is similar).
The natural operation among knots is the connected sum which is 
well-defined for the equivalence classes of knots. The sum operation is 
commutative and associative. Denote $K = K_1 \# K_2$ for the connected
sum of knots or the composite knot of $K_1$ and $K_2$. It turns out that
there is a nice algebraic topology method to compute the symplectic
Floer homology of the composite knot in terms of the ones of
$K_1$ and $K_2$. The method we used in this paper is in principle
the same for the instanton Floer homology \cite{li1, li3}.

It is easy to see that the special {\re}s of composite knot consists of two
different types: a single special {\re} arising from one knot and the unique
reducible special {\re} from the other, a circle of special {\re}s
arising from irreducible special {\re}s of both knots 
(see Proposition~\ref{ge}). The circle can also be interpreted as a gluing
parameters for two irreducible flat connections of the knot complement in
\cite{he}. We show that such a circle is a nondegenerate critical submanifold
of the symplectic action in the Bott sense \cite{bo}. A natural algebraic
topology method is to filter the critical submanifolds and to form a 
spectral sequence. This can not be done directly to our
symplectic Floer homology in \cite{li, li2} since the symplectic Floer
homology of knots is ${\Z}_{2N}$-graded in general ($N = N(K)$ is not
necessary zero). To overcome this difficulty, we define a $\Z$-graded
symplectic Floer homology of the braid which is an integral lifting of the one
in \cite{li} (see \cite{fs}). Using the monotonicity and some
properties of the special {\re} variety, we show that there is a well-defined
$\Z$-graded symplectic Floer homology of braids. The $\Z$-graded
symplectic Floer homology is only invariant under the Markov move of type
I and its inverse (Proposition~\ref{invariant}), not invariant under
the Markov move of type II and its inverse.
But there exists a nice relation between the $\Z$-graded 
and $\Z_{2N}$-graded symplectic Floer homologies. One of our main theorems
is the following.

\noindent{\bf Theorem A:} 
{\em 
\begin{enumerate}
\item There is a spectral sequence $(E^k_{n, j}(\p_{\b}), d^k)$ for $K = \ov{\b}$
with the $E^1_{n, j}(\p_{\b})$ term given by the $\Z$-graded 
symplectic Floer homology of the braid representing the knot $K$.
\item The spectral sequence $(E^k_{n, j}(\p_{\b}), d^k)$ converges to the
$\Z_{2N}$-graded symplectic Floer homology of the knot $K = \ov{\b}$.
\end{enumerate} }

Then we formulate a well-defined filtration for the
$\Z$-graded symplectic Floer chain complex of the braid representing
the composite knot
via the critical submanifolds. The filtration for the $\Z$-graded 
chain complex of the composite knot derives another spectral sequence.
Such a spectral sequence converges to the $\Z$-graded
symplectic Floer homology of the composite knot.
Using Theorem A, we obtain the $\Z_{2N}$-graded
symplectic Floer homology of the composite knot.

\noindent{\bf Theorem B:} 
{\em 
\begin{enumerate}
\item There is a spectral sequence 
$(E^r_{p,q}(\p_{\b}), d_r)$ determined by the filtration
(\ref{41}) for the braid $\b = \b_1 \Si^{n-1}(\b_2)$ of the composite knot
$K_1 \# K_2$.
\item The spectral sequence $(E^r_{p,q}, d_r)$ collapses at the
third term. Thus $E^3_{*,*} = E^{\infty }_{*,*}$
gives the ${\Z}$-graded Floer homology
$I^{(r_1, r_2)}_*(C_*^{(r_1, r_2)}(\b_1 \Si^{n-1}(\b_2)))$ of the braid $\b$.
\end{enumerate} }

By Theorem B, there are only two differentials needed to identify in order
to compute the $\Z$-graded symplectic Floer homology of the braid
representing the composite knot.
The differential $d_1$ basically consists of the integral differential 
on each braid $\b_j (j=1,2)$ and
plus two special boundary maps contributed from/to the unique
reducible special {\re} of $K_j (j=1,2)$. This is done by a cobordism
argument for the moduli spaces of $J$-holomorphic curves with variations
of Hamiltonian functions. We do not know how to characterize the differential
$d_2$ at this moment.

\noindent{\bf Theorem C:}
{\em The differential $d_1$ of the spectral sequence
$(E^r_{p,q}(\p_{\b}), d_r)$ in Theorem B is given by
 \[ d_1 = \partial_1^{(r_1)} \star \text{Id}_2 \pm \text{Id}_1 \star
\bd_2^{(r_2)}\pm
 d_{K_1} \star \text{Id}_2 \pm \text{Id}_1 \star d_{K_2}
+ \d_{K_1} \star \text{Id}_2 \pm \text{Id}_1 \star \d_{K_2}.\]
} 

See Definition~\ref{sbd} for the two special boundary maps
$d_{K_j}$ and $\d_{K_j}$ ($j=1, 2$).

This paper is organized as follows.
In \S 2.1, we briefly review the symplectic Floer homology of knots, 
and extend the ${\Z}_{2N}$-graded theory of knots
to the $\Z$-graded theory of braids in \S 2.2.
Using a spectral sequence, we link these two theories together in \S 2.3.
Theorem A is proved in Theorem~\ref{invariant}, Theorem~\ref{E1} and
Theorem~\ref{vds}. The filtration for the $\Z$-graded chain complex of the
braid representing the composite knot is 
formulated in \S 3.1. Theorem B is proved in \S 3.2 as
Theorem~\ref{ssi} and Proposition~\ref{44}. \S 4 devotes to the proof
of Theorem C (Theorem~\ref{boundary}).

\section{The symplectic Floer homology of a braid\label{braid}}

\subsection{The Floer homology of braids}

We briefly recall our definition of the Floer homology of braids in this
subsection. See \cite{li} for more details.

For any knot $K = \ov{\b}$ with $\b \in B_n$ braid group, 
the space ${\CR}(S^2 \setminus K)^{[i]}$ can be identified with the space
of $2n$ matrices $X_1, \cdots, X_n$, $Y_1, \cdots, Y_n$ in $SU(2)$ satisfying
\begin{equation} \label{zero}
\mbox{tr} (X_i) = \mbox{tr} (Y_i) = 0, \ \ \ \ \mbox{for $i =1, \cdots, n$},
\end{equation}
\begin{equation} \label{product}
X_1 \cdot X_2 \cdots X_n = Y_1 \cdot Y_2 \cdots Y_n .
\end{equation}
Let ${\CR}^*(S^2 \setminus K)^{[i]}$ be the subset of 
${\CR}(S^2 \setminus K)^{[i]}$ consisting of irreducible {\re}s. 
Note that $\pi_1(S^2 \setminus K)$ is generated by
$m_{x_i}, m_{y_i} (i =1, 2, \cdots, n)$
with one relation $\prod^n_{i=1} m_{x_i} = \prod^n_{i=1} m_{y_i}$.
There is a unique reducible
conjugacy class of {\re}s $s_K: \pi_1(S^3 \setminus K) \to U(1)$ is the
diagonal matrix
\[ s_K([m_{x_i}])  = s_K([m_{y_i}]) =  \begin{pmatrix}
i & 0 \\
0 & - i \\
\end{pmatrix}.\]

The space ${\CR}^*(S^2 \setminus K)^{[i]}$ is a monotone symplectic
manifold of dimension $4n-6$ by Lemma 2.3 in \cite{li}. 
The symplectic manifold $(M, \o)$ is called {\em monotone} if
$\pi_2(M) = 0$ or if there exists
a nonnegative $\a \geq 0$ such that
$I_{\o } = \a I_{c_1}$ on $\pi_2(M)$, where  
$I_{\o }(u) = \int_{S^2} u^*(\o ) \in {\R}$ and
$I_{c_1}(u) = \int_{S^2} u^*(c_1) \in \Z$ for $u \in \pi_2(M)$.
The braid $\b$ induces a diffeomorphism 
$\p_{\b }:  {\CR}^*(S^2 \setminus K)^{[i]} \to
{\CR}^*(S^2 \setminus K)^{[i]}$. The induced diffeomorphism $\p_{\b }$ is 
symplectic, and the fixed point
set of $\p_{\b }$ is ${\CR}^*(S^3 \setminus K)^{[i]}$ 
(see Lemma 2.4 in \cite{li}).

Let $H: {\CR}^*(S^2 \setminus K)^{[i]} \x {\R} \to {\R}$
be a $C^{\infty}$ time-dependent Hamiltonian function with 
$H(x, s) = H(\p_{\b}(x), s+1)$. Let $X_s$ be the corresponding vector field
from $\o (X_s, \cdot) = dH_s( \cdot, s)$, and 
\[\frac{d \psi_s}{ds} = X_s \circ \psi_s, \ \ \
\psi_0 = id.\]
Then we have $\psi_{s+1} \circ \p_{\b}^H = \p_{\b} \circ \psi_s$, where
$\p_{\b}^H = \psi_1^{-1} \circ \p_{\b}$.
Let $\O_{\p_{\b}}$ be the space of smooth contractible paths $\a$ in
${\CR}^*(S^2 \setminus K)^{[i]}$ such that $\a (s+1) = \p_{\b}(\a (s))$. 
The symplectic action 
$a_H: \O_{\p_{\b}} \to {\R}/\a 2N{\Z}$ is given by
\[ da_H(\g) \xi = \int_0^1 \o (\dot{\g} - X_s(\g), \xi) ds. \]
So the critical points of $a_H$ are the fixed points of $\p_{\b}^H$.
For $x \in \text{Fix}(\p_{\b}^H)$, define $\mu (x) = \mu_u (x, s_K) \pmod {2N}$,
where $\mu_u$ is the Maslov index and $N=N(K)$ is the first Chern number of 
the tangent bundle of ${\CR}^*(S^2 \setminus K)^{[i]}$.
The symplectic
Floer chain complex is defined by
\[C_j = \{ x \in \mbox{Fix}({\p}_{\b}) \cap {\CR}^*(S^2 \setminus K)^{[i]} :
\mu (x) = j \in {\Z}_{2N} \} . \] 
The following is the Proposition 4.1 and Theorem 4.2 in \cite{li}.
\begin{thm} \label{main}
For a knot $K= \ov{\b}$,
there is a well-defined ${\Z}_{2N}$-graded 
symplectic Floer homology $HF_*^{\mbox{sym}}(\p_{\b}^H)$.
The symplectic Floer homology 
$\{HF_j^{\mbox{sym}}(\p_{\b})\}_{j \in {\Z}_{2N}}$ is a knot invariant
and its Euler number is half of the signature of the knot $K$.
\end{thm}

\subsection{The $\Z$-graded symplectic Floer homology of braids}

In this subsection, we extend our ${\Z}_{2N}$-graded Floer homology
to the $\Z$-graded symplectic Floer homology of braids provided $\a N \neq 0$.
If $\a N = 0$, then the symplectic Floer homology of braids defined in
\cite{li} is $\Z$-graded.
Using a compatible filtration, we obtain a spectral sequence. 
The $E^1$ term of the spectral sequence
is the $\Z$-graded Floer homology of a braid,
and $E^{\infty}$ gives the
${\Z}_{2N}$-graded Floer homology.

\begin{lm} \label{pathc}
The symplectic manifold ${\CR}^*(S^2 \setminus K)^{[i]}$ is path connected.
\end{lm}
Proof: We prove this result by induction on the braid group $B_n$.
For $n=2$ (i.e., the knot $K$ can be represented by $K=\ov{\s_1^{\ve (K)}}$),
we see that ${\CR}^*(S^2 \setminus K)^{[i]}$ is a pillow case in \cite{li, lin}.
So the manifold ${\CR}^*(S^2 \setminus K)^{[i]}$ is path connected.
In particular, the space ${\CR}^*(S^2 \setminus K)^{[i]}$ is $K({\Z}^3, 1)$ 
space.

Suppose $n-1$ is true. For ${\CR}^*(S^2 \setminus K)^{[i]} = (H_n \setminus
S_n)/SU(2)$ in Lin's notation \cite{lin}, we have 
$(X_1, \cdots, X_n, Y_1, \cdots, Y_n)$ satisfies
\[ \text{tr} (X_j) = \text{tr} (Y_j) = 0, \ \  j=1, \cdots, n,\]
\begin{equation} \label{npr}
 X_1 \cdots X_n = Y_1 \cdots Y_n . \end{equation}
Applying the conjugate operation on $X_n$ and $Y_n$, we may assume that
\[ X_n = \begin{pmatrix} 
i & 0\\
0 & -i \\
\end{pmatrix}, \ \ \ 
Y_n = \begin{pmatrix}
i \cos \th  & \sin \th \\
- \sin \th  & - i \cos \th \\
\end{pmatrix}, \ \ \ 0 \leq \th \leq \pi .\]
If $\th = 0$ and $\pi$, then we get two copies of
$(H_{n-1} \setminus S_{n-1})/SU(2)$ which is path connected by the 
inductive hypothesis. For $0 < \th < \pi$, the equation (\ref{npr}) becomes
\[X_1 \cdots X_{n-1} 
\begin{pmatrix} 
\cos \th & - i \sin \th \\
- i \sin \th & \cos \th \\
\end{pmatrix} = Y_1 \cdots Y_{n-1}. \]
Let $R_{\th}$ be the representations in ${\CR}^*(S^2 \setminus K)^{[i]}$ 
satisfying the above equation. So the space $R_{\th }$ is non-singular
piece in ${\CR}^*(S^2 \setminus K)^{[i]}$. For $0 < \th , \th^{'} < \pi$, 
the space $R_{\th }$ is diffeomorphic to the space $R_{\th^{'}}$.
In particular, they are all diffeomorphic to $R_{\pi/2}$.
In this case, we see that ${\CR}^*(S^2 \setminus K)^{[i]}$ is a generalized
pillow case: 
\[{\CR}^*(S^2 \setminus K)^{[i]} = \bigcup_{0 \leq \th \leq \pi} R_{\th} .\]
It is clear that elements in $R_{\th}$ can be path connected to 
elements in $R_{\th^{'}}$.
So the manifold ${\CR}^*(S^2 \setminus K)^{[i]}$ is path connected.
\qed

I owe this idea to X-S. Lin.
For the path connected manifold ${\CR}^*(S^2 \setminus K)^{[i]}$ 
and $\a_i \in {\O }_{\p_{\b}}$ ($i=0, 1$),
there is a path $\g$ in ${\CR}^*(S^2 \setminus K)^{[i]}$ 
connecting $\a_0(0)$ to $\a_1(0)$. 
The space ${\O }_{\p_{\b}}$ is path connected so that 
$\pi_1({\O }_{\p_{\b}}, \a_0)$ 
is independent of a based point $\a_0 \in {\O }_{\p_{\b}}$.
Let $z_0 \in \text{Fix}\, (\p_{\b}) \subset {\O }_{\p_{\b}}$ be the
base point.
\begin{lm} \label{univ}
There exists a universal covering
 space ${\O}_{\p_{\b}}^*$ of ${\O }_{\p_{\b}}$.
\end{lm}
Proof:  The function space 
${\O }_{\p_{\b}}= \text{Map}(I, 0, 1; {\CR}^*(S^2 \setminus K)^{[i]}, z_0)$
has the homotopy type of a CW complex
and so an associated universal covering space. By Milnor's theorem 3.1
in \cite{mi}, there is a universal covering space
${\O}_{\p_{\b}}^*$ of ${\O }_{\p_{\b}}$. Note that over this CW complex there
is a universal covering space and so we can pull-back via the homotopy
equivalence to the covering space over 
$\text{Map}(I, 0, 1; {\CR}^*(S^2 \setminus K)^{[i]}, z_0)$.
As long as we work in the CW category this covering space is as good
as the universal covering space.
\qed

Using the standard algebraic topology in \cite{sp}, we can identify the
transformation group of the universal covering space of ${\O }_{\p_{\b}}$.
 
\begin{lm} \label{fund}
The transformation group of 
$\O_{\p_{\b}}^*$ is $\pi_2({\CR}^*(S^2 \setminus K)^{[i]}, z_0)$.
\end{lm}
Proof: Let $u(t) \in {\O }_{\p_{\b}}$ be a map from $I$ to 
${\CR}^*(S^2 \setminus K)^{[i]}$ with
$u(t)(0) = u(t)(1) = z_0$ and $u(t) (s)$ is contractible
path as a path in the variable $s$. Any loop $u(t) \in {\O }_{\p_{\b}}$ with
$u(0) = u(1) = z_0 \in \text{Fix}\, (\p_{\b})$ is a map
\begin{equation} \label{pi2}
u: (I \x I, I \x \{0, 1\} \cup \{0, 1\} \x I) \to 
({\CR}^*(S^2 \setminus K)^{[i]}, z_0).
\end{equation}
Note that the set of homotopy classes of maps in (\ref{pi2}) is in one-to-one
correspondence with 
$\pi_2({\CR}^*(S^2 \setminus K)^{[i]}, z_0)$ (see \cite{sp}).
\qed
 
The transformation group $\pi_2({\CR}^*(S^2 \setminus K)^{[i]}, z_0)$ 
is abelian. From the symplectic action map
$a_H: {\O}_{\p_{\b}} \to {\R}/\a (2N) {\Z}$,
we pullback the universal covering space of
${\R}/\a (2N) {\Z} \cong S^1$ over ${\O}_{\p_{\b}}$:
\[\begin{array}{ccc}
\ti{\O}_{\p_{\b}}=a_H^*({\R}) & \stackrel{\ti{a}_H}{\longrightarrow} & {\R}\\
\Big\downarrow & & \Big\downarrow \\
{\O}_{\p_{\b}}            & \stackrel{a_H}{\longrightarrow} &{\R}/\a (2N) {\Z}.
\end{array} \]
The pullback $a_H^*({\R})=\ti{\O}_{\p_{\b}} \to {\O}_{\p}$ 
is an infinite cyclic
subcovering space of the abelian universal covering space
${\O}_{\p_{\b}}^*$ with 
\[\pi_1(\ti{\O}_{\p_{\b}}) = {\Z} \triangleleft \pi_1({\O}_{\p_{\b}}^*) =
\pi_2({\CR}^*(S^2 \setminus K)^{[i]}, z_0). \]

Now the closed 1-form $da_H(z)$ has a functional $\ti{a}_H$
up to a constant such that
$\ti{a}_H : \ti{\O}_{\p_{\b}} \to \R$ is well-defined.
By adding a constant, we assume $\ti{a}_H (z_0) = 0$.
For a transformation element $g \in \pi_1(\ti{\O}_{\p_{\b}}) \triangleleft 
\pi_2({\CR}^*(S^2 \setminus K)^{[i]})$, we obtain 
\begin{equation} \label{deg}
\ti{a}_H (g(x)) = \ti{a}_H (x)
+ \deg (g) 2 \a N,
\end{equation}
where $\deg (g)$ is defined as $I_{\o}(g) = \deg (g) \a (2N)$.
Let $Im\, \ti{a}_H\, (\text{Fix}\, {\p_{\b}})$ be the image of
$\ti{a}_H$ of $\widetilde{\text{Fix}\, {\p_{\b}}}$; modulo $2 \a N \Z$,
$Im\, \ti{a}_H (\widetilde{\text{Fix}\, {\p_{\b}}})$ is a finite set. Thus
a set ${\R }_{\p_{\b}} = {\R } \setminus Im\, \ti{a}_H 
(\widetilde{\text{Fix}\, {\p_{\b}}})$
consists of the regular
values of the symplectic action $\ti{a}_H$ on $\ti{\O }_{\p_{\b}}$.
We
construct a $\Z$-graded symplectic Floer cohomology for
every $r \in {\R }_{\p_{\b}}$.

Given $x \in \text{Fix}\, {\p_{\b}} \subset \O_{\p_{\b}}$, let $x^{(r)} \in
\widetilde{\text{Fix}\, {\p_{\b}}} \subset
\ti{\O}_{\p_{\b}}$ be the unique lift of $x$ such that 
$\ti{a}_H\, (x^{(r)}) \in (r, r + 2\a N)$.
Note that if $\a N = 0$ then one already has an
$\Z$-graded symplectic
Floer homology as in \cite{li}, so we work on the case $\a N \neq 0$.
To cover the case $\a N \neq 0$, we may choose $r < \min \{a_H\,
(\text{Fix}\, \p_{\b}), \pmod {2 \a N}\}$ 
in the interval $(0, 2\a N)$. Then
there exists a unique element 
$x$ such that $\ti{a}_H (x^{(r)}) \in (r, r+ 2\a N)$.
Let $\mu^{(r)}(x) = \mu(x^{(r)}, s_K) \in \Z$ (the unique reducible conjugacy 
class $s_K$ of representations in ${\CR}(S^3 \setminus K)^{[i]}$) and define
the $\Z$-graded symplectic Floer chain group
\begin{equation} \label{cohain}
 C^{(r)}_n(\p_{\b}) = {\Z } \{ x \in \text{Fix}\, {\p_{\b}} \ \ \ | \ \ \
\mu^{(r)}(x) = n \in {\Z  }\} ,
\end{equation}
as the free $\Z$ module generated by $x \in \text{Fix}\, {\p_{\b}}$ 
with the lift
$x^{(r)}$ and $\mu(x^{(r)}, s_K) = n$.

If $\ov{z_0}$ ($\overline{s}_K$) is another choice of a based point $z_0$ 
(respectively $s_K$ for the Maslov index) and
$g(z_0) = \ov{z_0}$ ($g(s_K) = \overline{s}_K$) 
for some covering transformation $g \in \pi_1(\ti{\O}_{\p_{\b}})$, then the
corresponding choice of lift $\overline{x}^{(r)}$ of $x$ is just $g(x^{(r)})$
by the uniqueness.
Note that the integral Maslov index $\mu^{(r)}_u(x)$ is independent of the
choice of the based point $z_0$ ($s_K$) 
used in the definition of $a_H$ by (\ref{deg}).
We may choose $z_0 = s_K$ for simplicity on both lifts of the symplectic
action and the Maslov index (see \cite{fs}).
The following
lemma shows that the lift of the functional $a_H$
is compatible with a universal
lift of ${\R }/ 2 \a NZ (\cong S^1)$ via the monotonicity of the symplectic
manifold ${\CR}^*(S^2 \setminus K)^{[i]}$.

\begin{lm} \label{compatible}
The lift of the symplectic action over $\tilde{\O}_{\p_{\b}}$ is
compatible with the one of the Maslov index: 
for $g \in \pi_2({\CR}^*(S^2 \setminus K)^{[i]}) = \pi_1 (\ti{\O}_{\p_{\b}})$
with $\deg (g) = n$ in the sense of (\ref{deg}),
\[ \ti{a}_H (g(s_K)) = n \a (2N)  \ \ \ \text{if and only if} \ \ \ \
\mu^{(r)}(g(s_K), s_K) = n (2N). \]
\end{lm}
Proof: The result follows from the definition of
$\tilde{\O}_{\p_{\b}}$ and the monotonicity of 
$({\CR}^*(S^2 \setminus K)^{[i]}, \o)$ (see Lemma 2.3 in \cite{li}). 
So the Maslov index also
provides an integral lifting of element $x \in \text{Fix}\, (\p_{\b})$.
\qed
 
\begin{df} \label{integ} The integral Floer boundary map
$\bd^{(r)}: C^{(r)}_{n+1}(\p_{\b}) \to C^{(r)}_{n}(\p_{\b})$
is defined by
\[ \bd^{(r)}x = \sum_{y \in C^{(r)}_{n}(\p_{\b})} \# \hM (x, y) \cdot
 y, \]
where ${\CM }(x,y)$ denote the union of the components of 1-dimensional moduli
space of $J$-holomorphic curves and ${\hM }(x,y) = {\CM }(x, y)/{\R }$ is a
zero dimensional moduli space modulo the time translation.
The number $\# {\hM }(x, y)$
counts the points with sign in \cite{fh, li}.
\end{df}
Note that the boundary map $\bd^{(r)}$ only counts part of the 
boundary map in Proposition 4.1 of \cite{li}. We are going to show that
$\bd^{(r)} \circ \bd^{(r)} = 0$.
The corresponding homology groups are the $\Z$-graded symplectic
Floer homology $I^{(r)}_*(\p_{\b}), * \in \Z$.
See \cite{lo, sz}
for the following.

\begin{pro} \label{bub}
Suppose that elements in $\text{Fix}\, (\p_{\b})$ are nondegenerate.
(i) If $u \in {\CP}(x, y)$ for $x, y \in \text{Fix}\, {\p_{\b}}$
and $\tilde{u}$ is any
lift of $u$, then $\mu_{\tilde{u}}= \mu^{(r)}(y) - \mu^{(r)}(x)$.
 
(ii) Then there is a dense subset
${\J }_*(\p) \subset {\J }_{reg} (\p_{\b})$ of $\J$ such that (1) the zero
dimensional component of ${\hM }(x, y)$ is compact and (2) the one dimensional
component of ${\hM }(x^{'}, y^{'})$ is compact up to the splitting of two
isolated trajectories for $J \in {\J }_*(\p_{\b})$.
\end{pro}
 
Proposition~\ref{bub} plays the key role in showing that $\bd \circ \bd = 0$
(see \cite{li} \S 4).  We follow the same argument in \cite{fs, li}
to show that
${\bd }^{(r)} \circ {\bd }^{(r)} = 0$. 
 
\begin{lm} \label{26}
Under the same hypothesis in Proposition~\ref{bub},
$\bd^{(r)} \circ \bd^{(r)} = 0$.
\end{lm}
Proof: If $x \in C^{(r)}_{n+1}(\p_{\b})$, then by definition
the coefficient of $z \in C^{(r)}_{n-1}(\p_{\b})$  in
$\bd^{(r)} \circ \bd^{(r)}(x)$ is
\begin{equation}
 \sum_{y \in C^{(r)}_{n}(\p_{\b})} \# {\hM }(x, y) \cdot \# {\hM }(y,z).
\end{equation}
By Proposition~\ref{bub}, the boundary of the 1-dimensional
manifold $\hM (x, z) = {\CM }(x, z)/{\R}$ corresponds to two isolated
trajectories
${\CM }(x, y) \x {\CM }(y, z)$. Each term $\# \hM (x, y) \cdot \# \hM(y,z)$ is
the number of the 2-cusp trajectory of $\hM (x, z)$ with
$y \in C^{(r)}_{n}(\p_{\b})$.
For any such $y$ there are $J$-holomorphic curves $u \in {\CM }(x, y)$ and
$v \in {\CM }(y,z)$.
The other end of the corresponding
component of  the 1-manifold $\hM (x,z)$ corresponds to the splitting
${\CM }(x, y^{'}) \x {\CM }(y{'},z)$ with $u^{'} \in {\CM }(x, y^{'})$ and
$v^{'} \in {\CM }(y{'},z)$. It is impossible for $y{'}$ to be the reducible
$s$ because the $U(1)$ symmetry group would add one more parameter to the
moduli space.
Then $\hM (x, z)$ has an one
parameter family of paths
from $x$ to $z$ with ends $u \# v$ and $u^{'} \# v^{'}$ for appropriate
grafting. If we lift $u$ to $\tilde{u} \in \tilde{\CM} (x^{(r)},
\tilde{y})$, then
\begin{equation} \label{rest}
1 = \mu_{\tilde{u}} = \mu^{(r)}(x) - \mu^{(r)}(\tilde{y}) =
(n+1) - \mu^{(r)}(\tilde{y}).
\end{equation}
So $\mu^{(r)}(\tilde{y}) = n$; and $\tilde{y} = y^{(r)}$ is the preferred
lift, thus we have $\tilde{u} \in \tilde{\CM}(x^{(r)},
y^{(r)})$. Similarly $\tilde{v} \in \tilde{\CM}(y^{(r)}, z^{(r)})$.
Since $u^{'} \# v^{'}$ is homotopic to $u \# v$ rel $(x^{(r)}, z^{(r)})$, the
lift $\tilde{u}^{'} \# \tilde{v}^{'}$ is also a path with ends
$(x^{(r)}, z^{(r)})$. The symplectic action $a_H$ is
non-decreasing along the gradient trajectory $\tilde{u}^{'}$
with $\tilde{u}^{'}(- \infty ) = x^{(r)}$ and $\tilde{u}^{'}(+ \infty)
= \tilde{y}^{'}$, we have
\begin{equation} \label{ineq}
r < a(x^{(r)}) \leq a(\tilde{y}^{'}) \leq a(z^{(r)}) < r + \a 2N.
\end{equation}
By the uniqueness, $\tilde{y}^{'} = (y^{'})^{(r)}$ and using (\ref{rest})
for $u^{'}$, we have $\mu^{(r)}((y^{'})^{(r)}) = \mu^{(r)}(x^{(r)}) +1 = n$;
so $(y^{'})^{(r)} \in C^{(r)}_{n}(\p_{\b})$. Thus the algebraic number of
two-trajectories connecting $x^{(r)}$ and $z^{(r)}$ with index 2 is zero
by the orientation discussed in \cite{fh}.
\qed
 
Now $(C^{(r)}_{n}(\p_{\b}), \bd^{(r)}_n)_{n\in \Z}$ is indeed a
$\Z$-graded symplectic Floer chain complex of the braid $\b \in B_n$.
We call its homology to be a $\Z$-graded
symplectic Floer homology, denoted by
\[I^{(r)}_*(\p_{\b}) = H^*(C^{(r)}_*(\p_{\b}), \bd^{(r)}), \ \ \
* \in \Z. \]
From the construction we have that
\begin{align} \label{shift}
& \text{(1) if $[r,s] \subset {\R }_{\p}$, then
$I^{(r)}_*(\p_{\b}) = I^{(s)}_*(\p_{\b})$;}\\
 & \text{(2) $I^{(r)}_{* + (2N)}
(\p_{\b}) = I^{(r+ \a (2N))}_*(\p_{\b})$,}
\end{align}
where $\a (2N) ( > 0)$ is the minimal
number in $Im I_{\o }|_{\pi_2 ({\CR}^*(S^2 \setminus K)^{[i]})}$.
The relation between $I^{(r)}_*(\p_{\b})$ and
$HF_*^{\text{sym}}(\p_{\b})$ will be
discussed in \S 2.3. In the following, we notice that the $\Z$-graded
symplectic Floer homology $I_*^{(r)}(\p_{\b})$ 
is not a knot invariant in general.

Let $\{(J^{\lam }, H^{\lam }\}_{\lam \in \R}$ be an 1-parameter family
that interpolates from $(J^0, H^0)$ to $(J^1, H^1)$ and is constant in
$\lam $ outside $[0,1]$.
Define the perturbed $J$-holomorphic curve equations
\begin{equation} \label{jl}
{\pj }_{\lam } u_{\lam } (s, t ) =
\frac{\bd u_{\lam }}{\bd t} + J^{\lam}_s(
\frac{\bd u_{\lam }}{\bd s} - X_s^{\lam}(u_{\lam})) = 0 ,
\end{equation}
with the moving conditions $u_{\lam }(s+1, t) = \p^{\lam }_1( u_{\lam }(s, t))$
and the asymptotic values
\begin{equation} \label{ll}
\lim_{t \to - \infty} u_{\lam}(s, t) = x_0 \in \text{Fix}\, (\p_{\b, H^0}), 
\ \ \ \lim_{t \to + \infty} u_{\lam}(s, t) = 
x_1 \in \text{Fix}\, (\p_{\b, H^1}).
\end{equation}
This directly generalizes the $J$-holomorphic curve equation in the case
of $(J^0, \p_{\b, H^0})$ and $(J^1, \p_{\b, H^1})$ corresponding to
$\{(J^{\lam }, H^{\lam }\}_{\lam \in \R}$.
We define
\[ C^{(r)}_{\p_{\b}} = \min \{ \ti{a}_H(x^{(r)}) -r, 2 \a N +r - 
\ti{a}_H(x^{(r)}) | x \in
\text{Fix}\, (\p_{\b}) \},  \]
where $H$ is the Hamiltonian function which makes all $\text{Fix}\, (\p_{\b})$
nondegenerate and isolated. So the number $C^{(r)}_{\p_{\b}}$ is not
dependent on $(J^{\lam }, H^{\lam })$.

For each $x \in \text{Fix}\, (\p_{\b})$,
there is an open neighborhood $U_x$ in ${\O}_{\p_{\b}}$
such that (1) $U_x$ is evenly covered in
$\ti{\O}_{\p_{\b}}$, (2) for each $z \in U_x$,
$|a_{H^{\lam}}(z) - a_{H^{\lam}}(x)| < C^{(r)}_{\p_{\b}}/8$.
There are finite subcover $\{U_{x_1},
\cdots, U_{x_k}\}$ of $\text{Fix}\, (\p_{\b})$, 
and by Gromov's compactness theorem
\cite{fl, lo, sz}
we have $\ve_1 > 0$ such that if $\|Da_{H^{\lam}}(z)\|_{L^3_1} < \ve_1$ then
$z \in \bigcup_{i=1}^kU_{x_i}$. Let $\ve = \min \{ \ve_1,
C^{(r)}_{\p_{\b}}/8\}$. We set a deformation $\{J^{\lam}, H^{\lam}\}$
satisfying the usual perturbation requirements in \cite{fl, lo}, and also
satisfying
\begin{equation} \label{admi}
(i) |H_t^{\lam}(z)| < \ve /2, \ \ \ \
(ii) \|\n_J H_t^{\lam}(z)\|_{L^3_1} < \ve /2,
\end{equation}
for all $z \in {\O}_{\p_{\b}}$. These deformation conditions can be achieved by
the density statement in \cite{fl, lo, sz}.

Let ${\CP }_{1, \ve /2}$
be the set of $\{J^{\lam}, H^{\lam}\}$
which satisfies these extra conditions (\ref{admi}).
The moduli space ${\CM }_{\lam } (x, y )$
of (\ref{jl}) and (\ref{ll}) has the same analytic properties as the
moduli space ${\CM }(x, y)$ except for the translational invariance.
 
\begin{df} \label{cochain}
For each $n$ and $\{J^{\lam}, H^{\lam}\} \in {\CP }_{1, \ve /2}$, let
\[ \p_{01}^n (x_0) = \sum_{x_1 \in C^{(r)}_n (\p_{\b}; J^1, H^1)} \#
{\CM }^0_{\lam } (x_0, x_1) \cdot x_1, \]
be a homomorphism 
$\p_{01}^n : C^{(r)}_n (\p_{\b}; J^0, H^0) \to C^{(r)}_n (\p_{\b}; J^1, H^1)$,
where $\# {\CM }^0_{\lam } (x_0, x_1)$ is the algebraic number of the moduli
space ${\CM }^0_{\lam } (x_0, x_1)$ with compatible
orientation given by \cite{fh}.
\end{df}
 
\begin{pro} \label{p210}
For any continuation $(J^{\lam }, H^{\lam }) \in {\CP }_{1, \ve /2}$ which
is regular at the ends,
 
(i) the homomorphism $\{\p_{01}^*\}_{* \in Z}$ is a cochain map, i.e.,
\[ \bd^{(r)}_{n,1} \circ \p_{01}^n = \p_{01}^n \circ \bd^{(r)}_{n,0},\]
for all $n \in Z$;

(ii) for another regular pair $(J^{\lam_1}, H^{\lam_1})$ connecting
$(J^1, H^1)$ to $(J^2, H^2)$ in ${\CP }_{1, \ve /2}$,
\[\p^n_{12} \circ \p_{01}^n = \p^n_{02}, \ \ \ \p^n_{00} = \text{id}.\] 
\end{pro}
Proof: The proof follows from the same argument in the theorem 4 of
\cite{li} and
the method of the proof of Lemma~\ref{26}. See \cite{fl, fs}.
\qed

By Proposition~\ref{p210}, the $\Z$-graded symplectic Floer homology
$\{I^{(r)}_n(\p_{\b})\}_{n \in {\Z}}$ is invariant under the continuation of
$\{J^{\lam}, H^{\lam}\} \in {\CP }_{1, \ve /2}$.

\begin{pro} \label{invariant} 
The $\Z$-graded symplectic Floer homology 
$\{I^{(r)}_n(\p_{\b})\}_{n \in \Z}$ is invariant under the Markov move of type
I and its inverse.
\end{pro}
Proof: Suppose that we have a Markov move of type I:
change $\b$ to $\xi^{-1} \b \xi$ for some $\xi \in B_n$.
The element $\xi$ induces a diffeomorphism $\xi: Q_n \to Q_n$, and a symplectic
diffeomorphism $\xi \x \xi : Q_n \x Q_n \to Q_n \x Q_n$ which commutes with 
the $SU(2)$ action (see \cite{li} \S 4).
Let $f_{\xi}: {\CR}^*(S^2 \setminus K)^{[i]} 
\to {\CR}^*(S^2 \setminus K)^{[i]}$
be the symplectic diffeomorphism induced from $\xi \x \xi$. We have
\[ \p_{\b} = f_{\xi} \circ \p_{\xi^{-1} \b \xi} \circ f_{\xi}^{-1}\]
by changing variables via $f_{\xi}$. So $\ti{\p}_{\b} = 
f_{\xi} \circ \ti{\p}_{\xi^{-1} \b \xi} \circ f_{\xi}^{-1}$ for
Hamiltonian perturbations. Thus $f_{\xi}$ induces a natural identification
between $C_n^{(r)}(\p_{\b})$ and $C_n^{(r)}(\p_{\xi^{-1} \b \xi})$ since
$\text{Fix}(\p_{\xi^{-1} \b \xi})$ 
can be identified with
$\text{Fix}(\p_{\b})$ under $f_{\xi}$. Using the argument in 
Lemma~\ref{26} and \cite{fl}, the induced map $f_{\xi}$ is a chain map, i.e.,
it commutes with the $\Z$-graded boundary maps. There is an inverse chain map
$f_{\xi^{-1}}$ of $f_{\xi}$, and an isomorphism
\[ I_n^{(r)}(\p_{\b}) \cong I_n^{(r)}(\p_{\xi^{-1} \b \xi}).\]
It is clear that the argument goes through for the inverse operation of
a Markov move of type I. \qed

For the Markov move of type II (change $\b$ to $\s_n \b \in B_{n+1}$), there
is an imbedding $g: Q_n \x Q_n \to Q_{n+1} \x Q_{n+1}$ which induces
an imbedding
\[\hat{g}: \hat{H}_n (= {\CR}^*(S^2 \setminus \ov{\b})^{[i]}) \to 
\hat{H}_{n+1} (= {\CR}^*(S^2 \setminus \ov{\s_n \b})^{[i]}).\]
See \cite{li} for more details. The image $\hat{g}(\p_{\b})$ is invariant
under the operation of $\s_n$. The symplectic diffeomorphism $\p_{\s_n \b}$
can be connected with $\hat{g}(\p_{\b})$ through a Hamiltonian isotopy
$\psi_t$ in \cite{li}. Such a Hamiltonian isotopy $\psi_t$ is constructed by a
rotation about an axis. The Hamiltonian flow induced by $\psi_t$ connects
two general points in $S^2$. The Hamiltonian function $H$ is $x_3$ or $gx_3$
for $g \in SO(3)$. The function $H$ does not satisfy (\ref{admi}) for small
$\ve$ ($\n H= (0, 0, 1)$).

In \cite{sm}, Smale proves that the rotation group $SO(3)$ is a strong 
deformation retract of the orientation preserving $C^{\infty}$-diffeomorphisms
of $S^2$ with $C^r$ topology. Since the 2-sphere is simply connected, so every
symplectic diffeomorphism of identity component is Hamiltonian.
Thus any Hamiltonian isotopy connecting $\hat{g}(\p_{\b})$ with
$\p_{\s_n \b}$ can be obtained by the rotation about certain axis. I.e., the
Hamiltonian function is not a small perturbation for the Markov move of type
II and its inverse. So the continuation $(J^{\lam}, H^{\lam})$ realizing 
the Markov move of type II is not in ${\CP}_{1, \ve /2}$ in general. 
It is interesting to investigate this large perturbation via the method 
developed in \cite{her}. The homomorphism
$\p_{01}^n$ in Definition~\ref{cochain} is not well-defined. So the 
$\Z$-graded symplectic Floer homology of braids is not invariant under the
Markov move of type II and its inverse.
Although the $\Z$-graded symplectic Floer homology of braids is not a knot
invariant, we formulate a spectral sequence which converges to a knot
invariant in the next subsection.
The $\Z$-graded symplectic Floer homology of a braid $\b \in B_n$
is a stepping stone
to compute the ${\Z}_{2N}$-graded symplectic Floer homology of the knot
$K = \ov{\b}$.

\subsection{The spectral sequence for the symplectic Floer homology}

In this subsection
we are going to show that the $\Z$-graded symplectic Floer
homology $I_*^{(r)}(\p_{\b})$ for $r \in {\R }_{\p_{\b}}$
and $* \in \Z$ determines
the ${\Z }_{2N}$-graded symplectic Floer homology
$HF^*_{\text{sym}}(\p_{\b}) \, (* \in {\Z }_{2N})$.
The way to link them together is to filter the $\Z$-graded symplectic
Floer chain complex. The filtration, by a standard method
in algebraic topology,
 formulates a spectral sequence
which converges to the ${\Z }_{2N}$-graded symplectic Floer homology
$HF_*^{\text{sym}}(\p_{\b})$ of the knot $K = \ov{\b}$.
\begin{df} \label{filter}
For $r \in {\R }_{\p}, j \in {\Z }_{2N}$ and
$n \equiv j \ \ (\ mod \ \ (2N))$, define the free $\Z$ modules
\[ F_n^{(r)} C_j(\p_{\b}) = \sum_{k \geq 0} C^{(r)}_{n + (2N) k}
(\p_{\b}), \]
which gives a natural decreasing filtration on $C_*(\p_{\b})\, (* \in
{\Z }_{2N})$.
\end{df}
There is a finite length decreasing filtration of $C_j(\p_{\b})\,
(j \in {\Z }_{2N})$:
\begin{equation} \label{dea}
\cdots F_{n+(2N)}^{(r)}C_j(\p_{\b}) \subset F_n^{(r)}
  C_j(\p_{\b}) \subset F_{n - (2N) }^{(r)}C_j(\p_{\b})
  \cdots \subset C_j(\p_{\b}) .
\end{equation}
\begin{equation} \label{col}
C_j(\p_{\b}) = \bigcup_{n \equiv j \pmod {2N}}
F_n^{(r)} C_j(\p_{\b}).
\end{equation}
Note that the symplectic action is non-decreasing
along the gradient trajectories,
it follows that the boundary map
$\bd: F_n^{(r)} C_j(\p_{\b}) \to F_{n-1}^{(r)} C_{j-1}(\p_{\b})$
preserves the
filtration in Definition~\ref{filter}. The ${\Z}_{2N}$-graded
symplectic Floer cochain complex $(C_j(\p_{\b}), \bd)_{j \in
{\Z }_{2N}}$ has a decreasing bounded
filtration $(F_n^{(r)}C_*(\p_{\b}), \bd)$,

\begin{equation} \label{pre}
\begin{array}{ccccc}
 & \downarrow & \downarrow &  & \downarrow \\
\cdots & F_{n+ 2N}^{(r)}C_j(\p_{\b}) &
\subset F_n^{(r)}C_j(\p_{\b}) &
\cdots & \subset C_j(\p_{\b}) \\
 & \downarrow \bd^{(r)} & \downarrow \bd^{(r)} &  & \downarrow \bd \\
\cdots & F_{n+ (2N) -1}^{(r)}C_{j-1}(\p_{\b}) &
\subset F_{n-1}^{(r)}C_{j-1}(\p_{\b})
& \cdots & \subset C_{j-1}(\p; P, J) \\
& \downarrow & \downarrow &  & \downarrow
\end{array} .
\end{equation}
\begin{lm} \label{uds}
(1) The homology of the vertical chain subcomplex
$F_n^{(r)} C_*(\p_{\b})$ in the filtration (\ref{pre}) is
$F_n^{(r)} I_j^{(r)}(\p_{\b}).$
 
(2) There is a natural bounded filtration for
$\{I_*^{(r)}(\p_{\b}) \}_{* \in Z}$ the integral-graded symplectic
Floer homology,
\[\cdots F_{n+ (2N)}^{(r)}HF_j(\p_{\b}) \subset F_n^{(r)}
HF_j(\p_{\b}) \subset F_{n-(2N)}^{(r)}HF_j(\p_{\b})
\cdots \subset I_j^{(r)}(\p_{\b}),\]
where $F_n^{(r)}HF_j(\p_{\b}) = \text{Im}\, ( 
F_{n}^{(r)}I_j^{(r)}(\p_{\b}) \to I_j^{(r)}(\p_{\b}))$.
\end{lm}
Proof: The results follows from Definition~\ref{filter} and standard results
in \cite{sp} Chapter 9 and \cite{fs}. \qed
 
\begin{thm} \label{E1}
(i) There is a spectral sequence $(E^k_{n, j}(\p_{\b}), d^k)$ with
\[ E^1_{n,j}(\p_{\b}) \cong I_n^{(r)}(\p_{\b}), \ \ \
n\equiv j \ \ (mod \ \ (2N)), \]
and the higher differential $d^k :
E^k_{n, j}(\p_{\b}) \to
E^k_{n + (2N) k -1, j -1}(\p_{\b})$, and
\[E^{\infty }_{n,j}(\p_{\b}) \cong F_n^{(r)} HF_j^{\text{sym}}(\p_{\b})/
F_{n + (2N)}^{(r)}HF_j^{\text{sym}}(\p_{\b}) . \]
(ii) The spectral sequence $(E^k_{n, j}(\p_{\b}), d^k)$ converges to
$HF_*^{\text{sym}}(\p_{\b})$ (the ${\Z}_{2N}$-graded
symplectic Floer homology) of the knot $K = \ov{\b}$.
\end{thm}
Proof: (i) Note that
\[ F_n^{(r)} C_j(\p_{\b}) /F_{n + (2N)}^{(r)} C_j(\p_{\b}) =
C_n ^{(r)}(\p_{\b}) . \]
It is standard from \cite{sp} that there exists a spectral sequence
$(E^k_{n, j}, d^k)$ with $E^1$ term given by the homology of $
F_n^{(r)} C_j(\p_{\b}) /F_{n + (2N)}^{(r)} C_j(\p_{\b})$, so we have
$E^1_{n,j}(\p_{\b}) \cong I_n^{(r)}(\p_{\b})$
and $E^{\infty }_{n,j}(\p_{\b})$ is isomorphic to the bigraded $Z$-module
associated to the filtration $F^{(r)}$ of the $\Z$-graded symplectic
Floer homology $I_n^{(r)}(\p_{\b})$.
 
(ii) Since
$\text{Fix}\, {\p_{\b}}$ is a finite set and non-degenerate,
so the filtration $F$ is bounded
and complete from (\ref{dea}) and (\ref{col}).
Thus the induced spectral sequence converges to the
${\Z}_{2N}$-graded symplectic Floer homology.
Note that the grading is unusual
(jumping by $2N$ in each step), we list the terms for
$Z^k_{*, *}(\p_{\b})$ and $E^k_{*,*}(\p_{\b})$.
\begin{align*}
Z^k_{n,j}(\p_{\b}) & = \{ x \in F_n^{(r)}C_j(\p_{\b}) | \bd x \in
F_{n-1 + (2N) k}^{(r)} C_{j-1}(\p_{\b}) \}, \\
Z^{\infty}_{n,j}(\p_{\b}) & =
\{ x \in F_n^{(r)}C_j(\p_{\b})  | \bd x = 0 \}, \\
E^k_{n,j}(\p_{\b}) & = Z^k_{n,j}(\p_{\b})/
\{ Z^{k+1}_{n + (2N),j}(\p_{\b}) +
\bd Z^{k-1}_{n+ (k-1) (2N) +1,j+1}(\p_{\b})\}, \\
E^{\infty}_{n,j}(\p_{\b}) &= Z^{\infty}_{n,j}(\p_{\b}) /
\{ Z^{\infty}_{n + (2N),j}(\p_{\b}) +
\bd Z^{\infty}_{n+ (k-1) (2N) +1,j+1}(\p_{\b})\} . \end{align*}
Thus the Floer boundary map $\bd$ induces the higher differential
\[ d^k: E^k_{n, j}(\p_{\b}) \to
E^k_{n + (2N) k -1, j -1}(\p_{\b}) . \] \qed
 
\begin{thm}  \label{vds}
(1) For any continuation $(J^{\lam }, H^{\lam }) \in {\CP }_{1, \ve /2}$
which is
regular at ends, there exists an isomorphism
\[E^1_{n,j}({\p_{\b}}; J^0, H^0) \cong E^1_{n,j}({\p_{\b}}; J^1, H^1) .\]
(2) For each $k \geq 1$, $E^k_{n, j}(\p_{\b})$ are the invariants under the
continuation $(J^{\lam }, H^{\lam }) \in {\CP }_{1, \ve /2}$.
\end{thm}
Proof: By Theorem~\ref{E1}, we have an isomorphism
$E^1_{n,j}({\p_{\b}}; J^0, H^0) \cong I_n^{(r)}({\p_{\b}}; J^0, H^0)$,
so there exists an isomorphism by Proposition~\ref{p210}:
$I_n^{(r)}({\p_{\b}}; J^0, H^0) \to I_n^{(r)}({\p_{\b}}; J^1, H^1)$
which respects the filtration and
induces an isomorphism on the $E^1$ term. Since $E^1_{n, j}(\p_{\b})$
is an invariant under the 
continuation $(J^{\lam }, H^{\lam }) \in {\CP }_{1, \ve /2}$
by Proposition~\ref{p210} and Theorem~\ref{E1}, 
so (2) follows from Theorem 1 in \cite{sp} page 468.
\qed
 
If $\a N = 0$,
we already have the $\Z$-graded symplectic Floer homology
$HF_n^{\text{sym}}(\p_{\b})$ for $K = \ov{\b}$.
All these new spectral sequences should contain more information on
$({\CR}^*(S^2 \setminus K)^{[i]}; \p_{\b})$, 
they are also finer than $HF_*^{\text{sym}}(\p_{\b}),
 * \in {\Z}_{2N}$ the symplectic homology defined in \cite{li}.
If $\a N \neq 0$, we only have that $\{E^k_{n, j}(\p_{\b})\}$ is
invariant under the Markov move of type I and its inverse by 
Proposition~\ref{invariant} and Theorem 1 in \cite{sp} page 468.
In general $\{E^k_{n, j}(\p_{\b})\}_{1 \leq k < k(\p_{\b})}$ is not an 
invariant under the Markov move of type II and its inverse, where
$k(\p_{\b})$ is the minimal $k$ for
$E^k_{n, j}(\p_{\b}) = E^{\infty}_{n, j}(\p_{\b})$. This phenomena corresponds
to the quantum effect by the Markov move of type II and its inverse in the
symplectic Floer theory. All the higher differentials $d^k$ in the spectral 
sequence count the one dimensional moduli spaces of $J$-holomorphic curves
with larger energy. So the Markov move of type II and its inverse
will induce a spectral sequence homomorphism (actually isomorphism)
for large $k$, in particular $E^{\infty}$-term is invariant under the
Markov move of type II and its inverse (see \cite{li}).
In \cite{li4, oh}, these filtered information has been
discussed for monotone symplectic manifolds and monotone Lagrangians.

\begin{cor}
For $j \in {\Z }_{2N}$,
\[ \sum_{k \in Z} I^{(r)}_{ j + (2N) k }(\p_{\b})
= HF_j^{\text{sym}}(\p_{\b}) \]
if and only if all the differentials $d^k$ in the spectral sequence
$(E^k_{n, j}(\p_{\b}), d^k)$ are trivial (i.e. $k(\p_{\b}) = 1$).
In particular, the $\Z$-graded symplectic Floer homology of the
braid $\b$ is a knot invariant if $k(\p_{\b}) = 1$.
\end{cor}
In general, we see that \[\sum_{k \in \Z} I^{(r)}_{ j + (2N) k }(\p_{\b})
\neq HF_j^{\text{sym}}(\p_{\b})\] for $j \in {\Z}_{2N}$ for $\p_{\b}$.
The $\Z$-graded
symplectic Floer homology $I_*^{(r)}(\p_{\b})$ of the braid $\b$
can be thought as an integer lift of the
symplectic Floer homology $HF_*^{\text{sym}}(\p_{\b})$ of the knot
$K = \ov{\b}$.

\begin{cor} \label{eul}
The Euler characteristic of the spectral sequence is given by
\[\chi (E^k_{*, *}(\p_{\b})) = \chi (I^{(r)}_*(\p_{\b})) = 
\chi (HF_*^{\text{sym}}(\p_{\b})) = \frac{1}{2} \text{sign}\, (K).\]
\end{cor}
Proof: Note that the Euler characteristic of the chain complex is same as
the Euler characteristic of its induced homology. So 
the Euler characteristic $\chi (E^k_{*, *}(\p_{\b}))$ is independent of $k$.
In particular, 
\begin{eqnarray*}
\chi (E^k_{*, *}(\p_{\b})) & = &
\chi (E^1_{*, *}(\p_{\b})) = \chi (I^{(r)}_*(\p_{\b}))\\
& = & \chi (E^{\infty}_{*, *}(\p_{\b})) = \chi (HF_*^{\text{sym}}(\p_{\b})) \\
& = & - \lam_{CL}(K).
\end{eqnarray*}
Hence the result follows from Corollary 4.3 in \cite{li}. \qed

\section{The symplectic Floer homology of composite knots\label{sum}}

\subsection{The filtration of the composite knots}
 
We first recall the construction of composite knots. Let $K_j (j=1, 2)$ be a
knot in $S^3$. There is a neighborhood $U$ of $p_j \in K_j$ in $S^3$ such
that the pair $(U, U \cap K_j)$ is topologically equivalent to the
canonical ball pair $(B^3_j, B^1_j)$. By removing the
resulting pairs $(B^3_j, B^1_j)$ from $(S^3, K_j)$ and sewing the
resulting pairs by a homeomorphism $f: (\bd B_2^3, \bd B_2^1) \to  
(\bd B_1^3, \bd B_1^1)$, we form the pair connected sum
$(S^3, K_1) \# (S^3, K_2) = (S^3, K_1) \cup_f (S^3, K_2)$. In the convention,
we write the connected sum as $K = K_1 \# K_2$ for the connected sum of knots 
or the composite knot of $K_1$ and $K_2$.

If the knot $K_j$ is represented by a closure of a braid element $\b_j$,
then the composite knot $K_1 \# K_2$ can be described as the following.
If $K_1 = \ov{\b_1}, \b_1 \in B_n$ and $K_2 = \ov{\b_2}, \b_2 \in B_m$, then
\[K_1 \# K_2 = \ov{\b_1 \Si^{n-1} (\b_2)}, \ \ \ \ \
\b_1 \Si^{n-1} (\b_2) \in B_{n+m}, \]
where $\Si$ is the shift map on the inductive limit of the $B_n$'s and
$\Si (\s_i) = \s_{i+1}$. For example, the composite knot of the
figure $8$ $K_1 = \ov{\s_1 \s_2^{-1} \s_1 \s_2^{-1}}$ and the trefoil
$K_2 = \ov{\s_1^3}$ is $\ov{\b}$, where $\b = \s_1 \s_2^{-1} \s_1 \s_2^{-1}
\Si^{3-1}(\s_1^3) = \s_1 \s_2^{-1} \s_1 \s_2^{-1} \s_3^3$ (see \cite{bi}).

Let $K_1$ and $K_2$ be oriented knots in $S^3$, and
let $K = K_1 \# K_2$ be their connected sum. Let
$\mu_j \in \pi_1(S^3 \setminus K_j) (j =1, 2)$ be the meridian.
By Seifert-Van Kampen's theorem, we have
\[\pi_1(S^3 \setminus K) = \pi_1(S^3 \setminus K_1) *_f \pi_1
(S^3 \setminus K_2), \]
where the amalgamation homomorphism $f: \la \mu_1 \ra \to 
\la \mu_2 \ra$ is given by $f(\mu_1) = \mu_2$.
Let $\mu = \mu_1= \mu_2 \in \pi_1(S^3 \setminus K)$ be the meridian of $K$.
With the identification of the meridians, we can identify 
${\CR}^*(S^3 \setminus K)^{[i]}$ in the following.

\begin{pro} \label{ge}
${\CR}^*(S^3 \setminus K)^{[i]}$ is a disjoint union of
\begin{align*}
{\CR}^*(S^3 \setminus K_1)^{[i]} \star s_{K_2} & = \{[\r_1 * s_{K_2}] \in
{\CR}^*(S^3 \setminus K)^{[i]} \colon
s_{K_2} \, \text{reducible}\},\\
s_{K_1} \star {\CR}^*(S^3 \setminus K_2)^{[i]} & = \{[s_{K_1} * \r_2] \in
{\CR}^*(S^3 \setminus K)^{[i]} \colon s_{K_1} \, \text{reducible}\},\\
{\CR}^*(S^3 \setminus K_1)^{[i]} \star {\CR}^*(S^3 \setminus K_2)^{[i]} & =
\{[\r_1 * \r_2] \in {\CR}^*(S^3 \setminus K)^{[i]}\}.
\end{align*}
where $s_{K_j} (j=1, 2)$ is the unique reducible representation in
${\CR}(S^3 \setminus K_j)^{[i]}$.
For each $\r_j \in {\CR}^*(S^3 \setminus K_j)^{[i]}$, 
$\r_1 \star s_{K_2}$ and $s_{K_1} \star \r_2$ are single points, and
$\r_1 \star \r_2$ represents a $U(1)$ of inequivalent representation classes.
\end{pro}
Proof: Clearly the three spaces are disjoint. 
If $\r \in {\CR}^*(S^3 \setminus K)^{[i]}$, then we define
$\r_j = \r|_{\pi_1(S^3 \setminus K_j)}$ ($j=1, 2$) and write 
$\r = \r_1 \star \r_2$ with $\r_1 (\mu_1) = \r_2 (\mu_2)$. 
Given $\r_j \in {\CR}^*(S^3 \setminus K_j)^{[i]}$ we can form
\[\r = \r_1 * \r_2 \  \ \ \text{if and only if} \ \ \ 
\r_1 (\mu_1) = \r_2(\mu_2).\]
If both $\r_1$ and $\r_2$ are reducible, the
representation $\r_1 \star \r_2$ must be reducible. Hence
for $\r \in {\CR}^*(S^3 \setminus K)^{[i]}$ the restriction $\r_j$ we defined
must have at least one irreducible. 
For $\r_1 \in {\CR}^*(S^3 \setminus K_1)^{[i]}$, there is a unique abelian 
representation $s_{K_2}: \pi_1 (S^3 \setminus K_2)$ 
with $\r_1(\mu_1) = s_{K_2}(\mu_2)$
and $\text{tr}\, (s_{K_2}) = 0$. Thus there is a unique way to form
$\r_1 \star s_{K_2}$ and $s_{K_1} \star \r_2$ up to conjugacy.

For any $([\r_1], [\r_2]) \in {\CR}^*(S^3 \setminus K_1)^{[i]} \x 
{\CR}^*(S^3 \setminus K_2)^{[i]}$, there exists $g \in SU(2)$ such that
$g \r_2(\mu_2) g^{-1} = \r_1 (\mu_1)$. So $\r_1 \star (g\r_2 g^{-1}) \in 
{\CR}^*(S^3 \setminus K)^{[i]}$. The set of all cosets $[g] \in SU(2)/{\pm 1}$
satisfying
$g \r_2(\mu_2) g^{-1} = \r_1 (\mu_1)$ is a coset of the $U(1)$ subgroup
containing $[\r_2(\mu_2)] \in SU(2)/{\pm 1}$. Therefore
$[ \r_1 \star \r_2]$ is parameterized by $U(1)$ of inequivalent classes
since $SU(2)/{\pm 1}$ acts freely on the irreducible representations.
\qed 

Again there is a unique equivalent class of reducible representations
in ${\CR}(S^3 \setminus K)^{[i]}$. Such a reducible representation $s_K$
restricts to $s_{K_j} \in {\CR}(S^3 \setminus K_j)^{[i]}$ ($j=1, 2$).
Denote $s_K$ by $s_{K_1} \star s_{K_2}$.
The above result has a natural extension to the representations with fixed 
trace (see \cite{clm, he, li2}). This is known to Klassen in \cite{kl}.

By the identification between $\text{Fix}\, (\p_{\b})$ and
${\CR}^*(S^3\setminus K)^{[i]}$, we see that
${\CR}^*(S^3\setminus K)^{[i]}$ are critical points of $a_H$ on
$\O_{\p_{\b}}$, where $\O_{\p_{\b}}$ is the space of contractible
path $\g : {\R} \to {\CR}^*(S^2\setminus K)^{[i]}$ with
$\g(s+1) = \p_{\b}(\g (s))$ and $\b = \b_1 \Si^{n-1} (\b_2)$.

We define a nondegenerate critical submanifold as follows.
\begin{df} \label{ncm} 
A manifold $M$ is a {\em nondegenerate critical submanifold} of
$a_H$ on
$\O_{\p_{\b}}$ if (i) $da_H(x) = 0$ for every $x \in M$;
(ii) $T_xM = \ker \text{Hess}(a_H)(x)$;
(iii) the Maslov index $\mu (x) \pmod {2N}$ is constant on $M$.
\end{df}

Comparing with \cite{bo}, we only have the relative index for each
critical point in $M$. In our case, the negative normal bundle of
the critical submanifold is infinite dimensional. The well-defined
relative index (iii) can be used to replace the constant Morse index for
the critical submanifold.
Condition (ii) says that the manifold $M$ is nondegenerate in the normal
direction.

\begin{pro} \label{bott}
The $U(1)$ manifold $[\r_1 \star \r_2]$ is a nondegenerate
critical submanifold of $a_H$ on $\O_{\p_{\b}}$ in the
normal direction.
\end{pro}
Proof:
For any element $\r_1 \star (g\r_2g^{-1}) \in {\CR}^*(S^3\setminus K)^{[i]}$,
we have $da_H(\r_1 \star (g\r_2g^{-1})) = 0$ for all $g \in SU(2)/{\pm 1}$
with $\r_1 (\mu_1) = (g\r_2g^{-1})(\mu_2)$.
I.e., every point in $[\r_1 \star \r_2] \cong U(1)$ is a critical point.
The Hessian of $a_H$ at $\r_g = \r_1 \star (g\r_2g^{-1})$ is a 
self-adjoint operator
\[A_{\r_g} = J \n_{\r_g}: L^2_1(\r_g^*T{\CR}^*(S^2\setminus K)^{[i]})
\to L^2(\r_g^*T{\CR}^*(S^2\setminus K)^{[i]}).\]
By the deformation theory, the space $T_{\r}{\CR}^*(S^2\setminus K)^{[i]}$
(the group cohomology of $\pi_1(S^2 \setminus K)$ with twisted coefficients)
can be identified with the $\r^{[i]}$-twisted group cohomology
$H^1(S^2 \setminus K, \r^{[i]})$
of $S^2 \setminus K = S^2 \setminus (S^2 \cap K)$ 
since $S^2 \setminus K$ is a $K(\pi, 1)$-space.
Note that on the smooth part ${\CR}^*(S^2\setminus K)^{[i]}$ the
Zariski tangent space is the tangent space 
$T_{\r}{\CR}^*(S^2\setminus K)^{[i]}$.
For $\a \in H^1(S^2 \setminus K, \r^{[i]})$, then $\a$ is a tangent
to a path in the representation variety since all the obstructions (cup product
and higher Massey products) vanish.
So the kernel of $\ker A_{\r_g}$ is
$\ker A_{\r_g} = B^1(S^2 \setminus K, \r_g^{[i]})$ the
coboundary of 1-cochains in the twisted group cohomology.
For any $g=e^{u_0}$, we define
\[\r_t = \r_1 \star e^{-tu_0 + O(t^2)}\r_2 e^{tu_0+ O(t^2)}.\]
This provides an adjoint orbit containing $\r_g$, and the cocycle
corresponding to $\r_t$ is 
\[u = Ad \r_g u_0 - u_0.\]
So $u$ is the coboundary $\d u_0$. For any path $g_t = e^{tv_0}$, $\r_t = 
g_t^{-1}(\r_g)g_t$, the cocycle corresponding to the $\r_t$ is
$Ad \r_g  v_0 - v_0$.
On the other hand, $B^1(S^2 \setminus K, \r_g^{[i]})$ is the tangent
space of the adjoint orbit through $\r_g$.
\[B^1(S^2 \setminus K, \r_g^{[i]}) = T_{\r_g}[\r_1 \star \r_2]. \]
Therefore we obtain that $\ker A_{\r_g} = T_{\r_g}[\r_1 \star \r_2]$. So
$A_{\r_g}$ is nondegenerate along the normal bundle of
$[\r_1 \star \r_2]$ in $\O_{\p_{\b}}$. Note that
$\r_1 \star (g\r_2g^{-1})$ is equivalent to $(g^{-1} \r_1 g) \star \r_2$.
Hence the proof is complete.
\qed

Since the symplectic manifold ${\CR}^*(S^2\setminus K)^{[i]}$ is monotone 
by Lemma 2.3 in \cite{li}, so we have 
\[\a = \frac{\a_1 N(K_1) + \a_2 N(K_2)}{N(K)},\]
where $N(K) = g.c.d (N(K_1), N(K_2))$ and $a_H: \O_{\p_{\b}} \to
{\R}/{\a 2N {\Z}}$ is well-defined.
Next we verify the additivity for the Maslov index.
\begin{pro} \label{masl}
For $r_j \in {\R}_{\p_{\b_j}}$ and $\r_j \in {\CR}(S^3 \setminus K_j)^{[i]}$
($j=1, 2$), we have
\[\mu_K(\r_{10}^{(r_1)}\star \r_{2}^{(r_2)}, \r_{11}^{(r_1)}
\star \r_{2}^{(r_2)})
= \mu_{K_1}(\r_{10}^{(r_1)}, \r_{11}^{(r_1)}).\]
\[\mu_K(\r_1^{(r_1)}\star \r_{20}^{(r_2)}, \r_1^{(r_1)}\star \r_{21}^{(r_2)})
= \mu_{K_2}(\r_{20}^{(r_2)}, \r_{21}^{(r_2)}).\]
\end{pro}
Proof: Let $u$ be a path in $\ti{\O}_{\p_{\b}}$: $u(t, s) \in 
{\CR}(S^2 \setminus K)^{[i]}, (t, s) \in I \x I$ such that
$u(t, s+1) = \p_{\b}(u(t,s))$ for each $t \in [0,1]$, and
$u(0, s) = \r_1^{(r_1)}\star \r_{20}^{(r_2)}, 
u(1, s) = \r_1^{(r_1)}\star \r_{21}^{(r_2)}$.
By the very definition of the Maslov index, we have
\[\mu_K(\r_1^{(r_1)}\star \r_{20}^{(r_2)}, \r_1^{(r_1)}\star \r_{21}^{(r_2)})
= \mu_u(\bd (I \x I)),\]
the evaluation of the Maslov class in ${\CR}(S^2 \setminus K)^{[i]}$ on the
boundary of $I \x I$. It is deformation invariance rel. end points
(see \cite{clm1, fl}). The Maslov index 
$\mu_K(\r_1^{(r_1)}\star \r_{20}^{(r_2)}, \r_1^{(r_1)}\star \r_{21}^{(r_2)})$ is
an integral lifting of the one in $\O_{\p_{\b}}$.

For each $t$, $u(t, \cdot)$ is a contractible path in
${\CR}(S^2 \setminus K)^{[i]}$: $\{t\} \x I \to
{\CR}(S^2 \setminus K)^{[i]}$. In particular, there is $u_j(t, \cdot)$
is a contractible path in ${\CR}(S^2 \setminus K_j)^{[i]}$
such that $u_1(t, \cdot) = g(t) u_2(t, \cdot) g(t)^{-1}$ along the
meridian and $u(t, \cdot) = u_1(t, \cdot) \star u_2(t, \cdot)$.
So $u_1(t,s)$ satisfies $u_1(t, s+1) = \p_{\b_1}(u_1(t, s))$,
$u_1(0,s) = u_1(1, s) = \r_1^{(r_1)}$, and
$u_1(t, \cdot)$ is a contractible path in $\ti{\CR}(S^2 \setminus K_1)^{[i]}$.
Note that there are homotopy equivalent loops corresponding to
$\mu_{K_1} (\r_1, s_{K_1}) + 2N(K_1)k$. But there is a unique homotopy 
loop in $\ti{\CR}(S^2 \setminus K_1)^{[i]}$ corresponding
to $\mu_{u_1}(\r_1^{(r_1)}, \r_1^{(r_1)}) = 0$.
I.e., $u_1(t, s)$ on the lifting space $\ti{\CR}(S^2 \setminus K_1)^{[i]}$
is contractible to a constant loop $\ti{u_1} = \r_1^{(r_1)}$ since
$\mu_{K_1}(\r_1^{(r_1)}, \r_1^{(r_1)}) = 0$ by Lemma~\ref{compatible}.
Over the space ${\CR}(S^2 \setminus K_1)^{[i]}$, we have infinitely
many inequivalent loops at $\r_1$ with $\mu(\r_1^{(r_1)}, \r_1^{(r_1)})
\equiv 0 \pmod {2N(K_1)}$. For the unique lifting $\r_1^{(r_1)}$, we
obtain only one homotopy class for $\mu_{u_1}(\r_1^{(r_1)}, \r_1^{(r_1)}) = 0$.
Let $U_1(t, \t)$ be the homotopy of $u_1$ and $\ti{u}_1$. Thus
$U_1(t, \t)(\mu_1)$ parameterizes $g(t) u_2(t, \cdot) g(t)^{-1}(\mu_2)$
to $\r_2^{(r_2)}(\mu_2)$. Since $S^1$ is path connected, there exists
$g(t, \t)$ such that $U(t, \t) = U_1(t, \t) \star u_2(t, \cdot)$ with
$U_1(t, \t)(\mu_1) = g(t, \t)u_2(t, \cdot)g(t, \t)^{-1}(\mu_2)$ provides
a homotopy from $u$ to $\ti{u}_1 \star u_2$ in the space 
${\CR}(S^2 \setminus K)^{[i]}$. By the invariance of Maslov index, we have
$\mu_{U(t, 0)} = \mu_{U(t, 1)}$, i.e., $\mu_u = \mu_{\ti{u}_1 \star u_2}$.
The pullback $(\ti{u}_1 \star u_2)^*T({\CR}^*(S^2 \setminus K)^{[i]})$ over
$\bd (I \x I)$ is the fiber
product of a trivial Lagrangian on ${\CR}^*(S^2 \setminus K_1)^{[i]}$ with
a loop of Lagrangians from ${\CR}^*(S^2 \setminus K_2)^{[i]}$ with the proper
lifting property. By the Maslov class $\mu$,
\begin{equation*}
\begin{split}
\mu_u(\bd (I \x I)) & = \mu_{\ti{u}_1 \star u_2}(\bd (I \x I))\\
 & = \mu_{u_2}(\bd (I \x I))\\
 & = \mu_{K_2}(\r_{20}^{(r_2)}, \r_{21}^{(r_2)}).\\
\end{split} \end{equation*}
The other identity follows from the same argument. \qed

Note that the monotone symplectic manifold ${\CR}^*(S^2 \setminus K)^{[i]}$ is
not a product of the manifolds ${\CR}^*(S^2 \setminus K_1)^{[i]}$ and
${\CR}^*(S^2 \setminus K_2)^{[i]}$. So one can not apply the catenation
property directly (for the catenation property of the Maslov index 
see \cite{clm1}). Note that the singularity of the respresentaion space has
codimension bigger than 2, our Maslov index discussion can be proved in
the space $H_n$ (see \cite{li, lin} for the notation). Over the space 
$H_n \x H_m$ one can apply the catenation property of the Maslov index.
 
\begin{cor} \label{maslc}
For $r_j \in {\R}_{\p_{\b_j}}$ and $\r_j \in {\CR}^*(S^3 \setminus K_j)^{[i]}$
($j=1, 2$), we have
\begin{equation*}
\begin{aligned}
\mu_K(\r_1^{(r_1)}\star s_{K_2}) & = \mu_{K_1}(\r_1^{(r_1)}),\\
\mu_K(s_{K_1} \star \r_{2}^{(r_2)}) & = \mu_{K_2}(\r_{2}^{(r_2)}), \\
\mu_K(\r_1^{(r_1)}\star \r_{2}^{(r_2)}) & = \mu_{K_1}(\r_1^{(r_1)}) +
\mu_{K_2}(\r_{2}^{(r_2)}).\\
\end{aligned} \end{equation*}
\end{cor}
Proof: Note that the Malsov index on ${\CR}^*(S^2\setminus K)^{[i]}$
is calculated with respect to the unique reducible representation
$s_{K_1} \star s_{K_2} = s_K$.
So we have
\begin{eqnarray*}
\mu_K(\r_1^{(r_1)}\star s_{K_2}) & = &
\mu_K(\r_1^{(r_1)}\star s_{K_2}, s_{K_1} \star s_{K_2} )\\
& = & \mu_{K_1}(\r_1^{(r_1)}, s_{K_1}) \\
& = & \mu_{K_1}(\r_1^{(r_1)}).
\end{eqnarray*}
Similarly we obtain the following.
\begin{eqnarray*}
\mu_K(\r_1^{(r_1)}\star \r_{2}^{(r_2)}) & = &
\mu_K(\r_1^{(r_1)}\star \r_{2}^{(r_2)}, s_{K_1} \star s_{K_2})\\
& = & \mu_K(\r_1^{(r_1)}\star \r_{2}^{(r_2)}, s_{K_1} \star \r_{2}^{(r_2)})
+ \mu_K(s_{K_1} \star \r_{2}^{(r_2)}, s_{K_1} \star s_{K_2})\\
& = & \mu_{K_1}(\r_1^{(r_1)}) + \mu_{K_2}(\r_{2}^{(r_2)}).\end{eqnarray*}
The second equality follows from the path additivity of the Maslov index
(see \cite{clm1}), and the 
last equality from Proposition~\ref{masl} and the previous two
identity. \qed

\begin{cor} \label{hess}
For each $g \in U(1)$ and $\r_j \in {\CR}(S^3 \setminus K_j)^{[i]} (j=1,2)$, 
we have
\[ \mu_K(\r_1 \star \r_2) \equiv 
\mu_K(\r_1 \star_g \r_2) \equiv \mu_{K_1}(\r_1) + \mu_{K_2}(\r_2) \pmod {2N(K)}.\]
\end{cor}
Proof: The proof follows the same argument from Proposition~\ref{masl} with
$\pmod {2N(K_j)}$.
So we get
\begin{equation*} 
\begin{split}
\mu_K(\r_1 \star_g \r_2)& \equiv \mu_{K_1}(\r_1^{(r_1)}) \pmod {2N(K_1)}
+ \mu_{K_2}(\r_2^{(r_2)}) \pmod {2N(K_2)}\\
& \equiv \mu_{K_1}(\r_1) + \mu_{K_2}(\r_2) \pmod {2N(K)}.\\
\end{split}  \end{equation*}
Hence it is independent of the element $g \in U(1)$.
\qed

By Proposition~\ref{bott} and Corollary~\ref{hess}, we have showed that
the $U(1)$-critical submanifold $[\r_1 \star \r_2]$ is a nondegenerate
critical submanifold in the sense of Definition~\ref{ncm} for
$\r_j \in {\CR}^*(S^3 \setminus K_j)^{[i]}$ ($j=1, 2$).
By Proposition~\ref{masl}, we see that $[\r_1^{(r_1)} \star \r_2^{(r_2)}]$
is also a nondegenerate critical submanifold in the sense of
Definition~\ref{ncm} for the $\Z$-graded symplectic Floer homology.

\begin{df} \label{icc}
The integral chain complex $C_n^{(r_1, r_2)}(\p_{\b_1 \Si^{n-1}(\b_2)})$ 
of the braid $\b_1 \Si^{n-1}(\b_2)$
with $K = K_1 \# K_2 = \ov{\b_1 \Si^{n-1}(\b_2)}$ is defined to be 
\[C_n^{(r_1, r_2)}(\p_{\b_1 \Si^{n-1}(\b_2)})=\{
\r_1\star \r_2 \in \text{Fix}\, (\p_{\b})= {\CR}^*(
S^3 \setminus K)^{[i]} |
\mu_K(\r_1^{(r_1)}\star \r_2^{(r_2)}) = n \]
\[\text{and}\, \ti{a}_{H_j}(\r_j^{(r_j)}) \in (r_j, r_j + 2 \a_j N(K_j),
r_j \in {\R}_{\p_{\b_j}}\}.\]
\end{df}

Note that $\text{Fix}\, (\p_{\b})$ $(\b = \b_1 \Si^{n-1}(\b_2))$ consists of
$\r_1 \star s_{K_2}, s_{K_1} \star \r_2$ and $\r_1 \star \r_2$. There is an
one-to-one correspondence between the elements in $\ti{\text{Fix}}\, (\p_{\b})$
and the elements in
\[\{\r_1^{(r_1)}\star \r_2^{(r_2)} | \r_j^{(r_j)} \in {\CR}(S^3 \setminus
K_j)^{[i]}, \r_1 \star \r_2 \in {\CR}^*(S^3 \setminus K)^{[i]},
\ti{a}_{H_j}(\r_j^{(r_j)}) \in (r_j, r_j + 2\a_j N(k_j))\}.\]
By Proposition~\ref{masl}, we have that Definition~\ref{icc} for
the composite knot is well-defined.
Fixing the symplectic action with respect to the 
$r_j \in {\R}_{\p_{\b_j}}$ is necessary.
Since there are $g_{j, k} \in \pi_1(\ti{\O}_{\p_{\b_j}})$ with
$\ti{a}_{H_j}(g_{j, k} \r_j^{(r_j)}) \in (r_j + 2\a_j k N(K_j),
r_j + 2\a_j (k+1)N(K_j))$ and $\mu_{K_j}(g_{j, k} \r_j^{(r_j)})=
\mu_{K_j}(\r_j^{(r_j)}) + 2k N(K_j)$, we have
\begin{equation*} \begin{split}
\mu_K(g_{1,kN_2} \r_1^{(r_1)}\star g_{2, -kN_1} \r_2^{(r_2)}) & = 
\sum_{j=1}^2 \mu_{K_j}(g_{j, (-1)^{j+1}kN_{-j+3}} \r_j^{(r_j)})\\
& = \sum_{j=1}^2 (\mu_{K_j}(\r_j^{(r_j)}) + 2(kN_2) N(K_1) + 2(-kN_1) N(K_2) \\
& = \sum_{j=1}^2 (\mu_{K_j}(\r_j^{(r_j)}).
\end{split} \end{equation*}
So $\mu_K(g_{1,kN_2} \r_1 \star g_{2, -kN_1} \r_2) \equiv 
\sum_{j=1}^2 \mu_{K_j}(\r_j^{(r_j)}) \pmod {2N(K)}$ for infinitely many $k$.
I.e., the element 
$g_{1,kN_2} \r_1^{(r_1)} \star g_{2, -kN_1} \r_2^{(r_2)}$ has the same Maslov
index $n$, but different symplectic actions with respect to $r_j (j=1, 2)$.
Thus our chain groups $C_*^{(r_1, r_2)}(\p_{\b})$ is well-defined, finitely
generated from $C_*^{(r_j)}(\p_{\b_j}) (j=1, 2)$, and
the integral liftings of $C_*(\p_{\b_1 \Si^{n-1}(\b_2)})$ corresponding to
$r_j \in {\R}_{\p_{\b_j}}$ ($j=1, 2$) is the integral
chain groups $C_*^{(r_1, r_2)}(\p_{\b_1 \Si^{n-1}(\b_2)})$.
See \cite{li3, li4} for more discussions on this kind of formulations.
\subsection{Spectral sequences for the composite knots representing by braids}

In this subsection, we form a filtered Floer chain complex for the 
braid $\b = \b_1 \Si^{n-1}(\b_2)$ representing the
composite knot $K = K_1 \# K_2$. The filtered chain complex of the braid
naturally associates to a spectral sequence. The spectral sequence
converges to the $\Z$-graded symplectic Floer homology
$I^{(r_1, r_2)}_*(\p_{\b})$ and collapses at third term.

For the composite knot $K_1 \# K_2 = \ov{\b_1 \Si^{n-1}(\b_2)}$ 
as in \S 3.1, we know that
the $\Z$-graded symplectic Floer chain groups $C^{(r_1, r_2)}_*(\p_{\b})$
are generated freely by
\[\r_1 \star s_{K_2}, \ \ \ \
s_{K_1} \star \r_2, \ \ \ \ \r_1 \star \r_2.\]
The free generated group
$C^{(r_1, r_2)}_*(\p_{\b}) = \oplus_q C_q^{(r_1, r_2)}(\p_{\b})$ is
a graded differential group with respect to the boundary map
$\bd^{(r_1, r_2)}_*$,
$\bd C_q^{(r_1, r_2)}(\p_{\b}) \subset
C_{q-1}^{(r_1, r_2)}(\p_{\b})$. Recall the
$\Z$-graded symplectic Floer chain group 
\[C_q^{(r_1, r_2)}(\p_{\b}) = \{ \r \in
{\CR}^*(S^3 \setminus K)^{[i]} \ \ | \mu (\ti{\r}) = q, \] \[ \ti{a}_{H_j}
({\r_j}^{(r_j)})
\in (r_j, r_j + 2 \a_j N_j) \ \text{for $\r = \r_1 \star \r_{K_2}$ and
$\r = \r_1 \star s_2$} \}. \]

There is an associated
filtration compatible with the grading, 
defined as follows:
\begin{equation} \label{41}
 F_pC^{(r_1, r_2)}_*(\p_{\b})= \bigoplus_{k \leq p} \{ \Z \la \r \ra \ \ | \ \
\r \in C^{(r_1, r_2)}_k(\p_{\b})\}.
\end{equation}
The filtration is an increasing filtration
$F_pC_q^{(r_1, r_2)}(\p_{\b}) \subset F_{p+1}C_q^{(r_1, r_2)}(\p_{\b})$.
The filtration (\ref{41}) is different
from the decreasing filtration (\ref{filter}).
Again the symplectic action functional is non-decreasing along
the gradient flows ($J$-holomorphic curves).
It follows that the $\Z$-graded symplectic Floer
boundary map $\bd^{(r_1,r_2)}: F_pC_q^{(r_1, r_2)}(\p_{\b})
\hookrightarrow F_pC_{q-1}^{(r_1, r_2)}(\p_{\b})$
preserves the filtration (\ref{41}). By the
standard method in algebraic topology \cite{sp}, we have a spectral sequence
induced from the filtration (\ref{41}).
 
\begin{lm} \label{generate}
For the braid $\b = \b_1 \Si^{n-1}(\b_2)$, we
have the $E^0_{*,*}$ and $E^1_{*, *}$ terms of the filtration (\ref{41})
given by the following.
\begin{equation*}
\begin{split}
E^0_{p,q} & =  F_pC_q^{(r_1, r_2)}(\p_{\b})/F_{p-1}C_q^{(r_1, r_2)}(\p_{\b})
= C_p^{(r_1, r_2)}(\p_{\b}) ; \\
E^1_{p,0} & =  \oplus Z \la \r \ra, \ \ \mbox{where} \ \
\r \in C^{(r_1, r_2)}_p(\p_{\b}), \ \ p =
\mu_K(\ti{\r});\\
E^1_{p,q} & =  0, \ \ \mbox{for} \ q\neq 0 \ \ \mbox{and} \ \
\r = \r_1 \star s_{K_2}, s_{K_1} \star \r_2; \\
E^1_{p,1} & =  \oplus Z \la (\r_1 \star \r_2)_1 \ra, \ \ \
\mu_K(\r_1^{(r_1)} \star \r_2^{(r_2)}) = p; \\
E^1_{p,q} & =  0, \ \ \ \mbox{for} \ \ q \neq 0, 1, \ \
\r = \r_1 \star \r_2, p = \mu_K (\r_1^{(r_1)} \star \r_2^{(r_2)}),
\end{split}
\end{equation*}
where $(\r_1 \star \r_2)_1$ is the critical point of the Morse index 1 in 
$U(1) = \r_1 \star \r_2$.
\end{lm}
Proof: The $E^0_{*, *}$-terms follow from the very definition of spectral
sequence induced from the filtration (\ref{41}).
The $E^1_{*,*}$-terms follows from all three types of
generators and $H_*(U(1), \Z)$, since $d_0$ is
the differential of $\r$'s for $\r \in C^{(r_1, r_2)}_p(\p_{\b})$.
If $\r = \r_1 \star s_{K_2}$ and $s_{K_1} \star \r_2$,
it is a single point so that
$E^1_{p, 0}= Z\la \r \ra$ for $p = \mu(\ti{\r})$ and zero otherwise.
If $\r = \r_1 \star \r_2$,
then it gives an $U(1)$ component which
$d_0$ is the boundary map of
the standard cellular chain complex of $U(1)$.
So the $E^1_{*, *}$-terms follows from the Maslov index calculations
and the homology of $U(1)$ for each copy of $\r_1 \star \r_2$.
\qed

One can define a precise perturbation of the symplectic action such that
all the critical points are isolated (see \cite{bo, li3}). The perturbation 
keeps the isolated points $\r_1 \star s_{K_2}$ and $s_{K_1} \star \r_2$,
and generates a Morse function along the tubular neighborhood of the
critical submanifold $\r_1 \star \r_2$. This also gives an interpretation
of the $E^1_{*, *}(\p_{\b})$ term in Lemma~\ref{generate}.

\begin{thm} \label{ssi}
There is a spectral sequence $(E^r_{p,q}(\p_{\b}), d_r)$ 
determined by the filtration
(\ref{41}) for the braid $\b = \b_1 \Si^{n-1}(\b_2)$ of the 
composite knot $K_1 \# K_2$.
The spectral sequence with
($E^1_{p, q}$ given by Lemma~\ref{generate})
\begin{equation}
 E^1_{p,q}(\p_{\b}) \cong H_{p+q}(F_pC^{(r_1, r_2)}_*(\p_{\b})/F_{p-1}
C^{(r_1, r_2)}_*(\p_{\b})),
\end{equation}
converges to the $\Z$-graded symplectic Floer homology
$I^{(r_1, r_2)}_*(\p_{\b}) (* \in \Z)$ of the braid ${\b}$.
\end{thm}
Proof: Because ${\CR}^*(S^3 \setminus \ov{\b})^{[i]}$ is compact, 
so $E^1_{p,q}(\p_{\b})$ is
finitely generated by Lemma~\ref{generate}.
In particular the filtration is bounded. Since
$\bigcup_p F_p C^{(r_1, r_2)}_*(\p_{\b}) = C^{(r_1, r_2)}_*(\p_{\b})$
(the filtration is exhaustive),
we get a convergent spectral sequence by the result in \cite{sp}.
By the definition of the filtration $F_pC^{(r_1, r_2)}_*(\p_{\b})$,
the higher differentials calculate all the integral Floer boundary map
$\bd^{(r_1, r_2)}_{\p_{\b}}$ (see \cite{fs, li3, li4}). 
Thus the spectral sequence converges to the
$\Z$-graded symplectic Floer homology of the braid ${\b}$.
\qed

\begin{pro}\label{44}
The spectral sequence $(E^r_{*,*}(\p_{\b}), d_r)$ in Theorem~\ref{ssi}
collapses at the third term. Thus the term 
$E^3_{*,*}(\p_{\b}) = E^{\infty }_{*,*}(\p_{\b})$
gives the ${\Z}$-graded symplectic Floer homology
$I^{(r_1, r_2)}_*(C_*^{(r_1, r_2)}(\p_{\b}))$ of the braid 
$\b = \b_1 \Si^{n-1}(\b_2)$.
\end{pro}
Proof: For the differential 
$d_r: E^r_{p,q}(\p_{\b}) \to E^r_{p-r, q+r-1}(\p_{\b}), r\geq 3$,
we have $E_{p, q}^r(\p_{\b}) = 0$ for $q < 0$ or $q > 1$ by 
Lemma~\ref{generate}.
The differential $d_r$ has either source zero or target zero.
So $d_r = 0$ for all $r \geq 3$.
Hence the spectral sequence collapses at the third term.
\qed
 
By Theorem~\ref{E1}, we can form another spectral sequence from $I^{(r_1,
r_2)}_*(\p_{\b})$ which converges to the ${\Z}_{2N}$-graded Floer homology
$HF_*^{\text{sym}}(K_1 \# K_2)$, where $N$ is the greatest common divisor of
$N(K_1)$ and $N(K_2)$.

In order to make some computations of $HF_*^{\text{sym}}(K_1 \# K_2)$,
we need to describe all the differentials in
Theorem~\ref{ssi} and combine the differentials in Theorem~\ref{E1}.
By Lemma~\ref{44}, we know that there are only two possibly nontrivial
differentials $d_1$ and $d_2$. In the section 4, we give a description of
$d_1$. We do not know how to characterize the differential $d_2$ yet. 
We leave it to a future study.

We define the associated Poincar\'{e}-Laurent polynomial for the spectral 
sequence by
\[P^{(r)}_K(E^k_{*,*},t) = \sum_{n \in \Z}(\dim_{\Z}E^k_{n, j})t^n. \]
By Theorem~\ref{E1}, $P^{(r)}_K(E^k_{*,*},t) = \sum_{n \in \Z}(\dim_{\Z}
I_n^{(r)}(\p_{\b})) t^n$.
\begin{pro} \label{twoid}
The following identities hold.
\begin{equation} \label{two1}
P^{(r)}_K(E^k_{*,*},t) = (1+t^{2Nk-1})P^{(r)}_K(B^k_{*,*},t) +
P^{(r)}_K(E^{k+1}_{*,*},t).
\end{equation}
\begin{equation} \label{two2}
P^{(r)}_K(E^1_{*,*},t) = \sum_{i=1}^k(1+t^{2Ni-1})P^{(r)}_K(B^i_{*,*},t) +
P^{(r)}_K(HF^{\text{sym}}_*,t).
\end{equation}
\end{pro}
Proof: Let $Z_{n,j}^k = \ker \{d^k: E^k_{n,j} \to E^k_{n+2Nk-1, j-1}\}$
and $B_{n,j}^k = \text{Im}\, d^k \cap E^k_{n,j}$. There are two short 
exact sequences:
\[0 \to Z_{n,j}^k \to E^k_{n,j} \to B^k_{n+2Nk-1, j-1} \to 0;\]
\[0\to B_{n,j}^k \to Z_{n,j}^k \to E^{k+1}_{n,j} \to 0.\]
So the degree $2Nk-1$ of $d^k$ derives the following identity.
\[P^{(r)}_K(E^k_{*,*},t) = (1+t^{2Nk-1})P^{(r)}_K(B^k_{*,*},t) +
P^{(r)}_K(E^{k+1}_{*,*},t).\]
Note that the spectral sequence converges at finite steps.
So the identity (\ref{two2}) follows from by iterating the identity
(\ref{two1}) and using Theorem~\ref{E1}. \qed

Note that $\c_K(E^k_{*,*}(\p_{\b})) = P^{(r)}_K(E^k_{*,*}(\p_{\b}), -1)$ and
Corollary~\ref{eul} is an easy consequence of Proposition~\ref{twoid}.
By Lemma~\ref{generate} and Theorem~\ref{ssi}, we obtain
\begin{equation} \label{add1}
\c_{K_1 \# K_2}(E^1_{*,*}(\p_{\b})) = \c_{K_1}(E^1_{*,*}(\p_{\b_1})) +
\c_{K_2}(E^1_{*,*}(\p_{\b_2})).
\end{equation}
By Corollary~\ref{eul} and (\ref{add1}),
\[\c_{K_1 \# K_2}(E^k_{*,*}(\p_{\b})) = \c_{K_1}(E^k_{*,*}(\p_{\b_1})) +
\c_{K_2}(E^k_{*,*}(\p_{\b_2})).\]
In particular, we have the equality on the $E^{\infty}$-term:
\[\c_{K_1 \# K_2}(HF^{\text{sym}}_*(\p_{\b})) =
\c_{K_1}(HF^{\text{sym}}_*(\p_{\b_1})) + 
\c_{K_2}(HF^{\text{sym}}_*(\p_{\b_2})).\]
Thus we obtain the classical result
\begin{equation}
\text{sign}(K_1\# K_2) = \text{sign}(K_1) + \text{sign}(K_2),\end{equation}
as one of identities from our spectral sequence for the symplectic
Floer homology of the braid.

\section{Description of $d_1$}

We are going to define two special boundary maps which are not used in the
definitions of the $\Z$-graded and the
${\Z}_{2N}$-graded symplectic Floer homologies in \cite{li} and \S 2 of the
present paper.
These two special boundary maps do contribute
for the composite knot. In particular, these two
special boundary maps are part of the differential
$d_1$ in the spectral sequence in Theorem~\ref{ssi} and Proposition~\ref{44}. 
In the gauge theory, we have
the same special boundary maps used in the instanton 
Floer homology of connected sums of two integral homology 3-spheres
\cite{li1, li3}.
 
\begin{df} \label{sbd}
$d_{\b} : C_1^{(r)}(\p_{\b}) \to \Z \la s_K \ra$ is defined by
\[d_{\b} \r = \# \hat{{\CM}}^1(\r, s_K) \cdot s_K, \]
and $\delta_{\b} : \Z \la s_K \ra \to C_{-1}^{(r)}(\p_{\b})$ is defined by
\[ \delta_Y s_K  = \sum_{\r \in C_{-1}(\p_{\b})} \#
\hat{{\CM}}^1(s_K, \r) \cdot \r. \]
\end{df}
 
Note that by fixing $\mu^{(r)}(s_K)=0$, both maps $d_{\b}$ and
$\d_{\b}$ are counting the one dimensional moduli spaces from/to
the unique reducible representation $s_K$.
 
\begin{lm} \label{special}
$ d_{\b} \partial^{(r)} = 0, \ \ \partial^{(r)} \delta_{\b} = 0$.
\end{lm}
{\bf Proof:} The proof is similar to the
proof that
$\bd^{(r)} \circ \bd^{(r)} = 0$ in Lemma~\ref{26} for the 
$\Z$-graded symplectic Floer chain complex of a braid
(also see \cite{fl, li}).  \qed
 
Note that for the Floer boundary map $\bd$ in \cite{li} we also have
$d_{\b} \circ \bd = 0$ and $\bd \circ \delta_{\b} = 0$.
Let $\bd_j^{(r_j)}: C_*^{(r_j)}(\p_{\b_j}) \to
C_{* -1}^{(r_j)}(\p_{\b_j}) (j =1, 2)$
be the boundary map in Definition~\ref{integ}
of the $\Z$-graded symplectic Floer chain complex of the braid $\b_j$.
The main result of this section is a description of $d_1$
in Theorem~\ref{ssi}.
 
\begin{thm} \label{boundary}
The differential $d_1$ of the spectral sequence
$(E^r_{p,q}(\p_{\b}), d_r)$ in Theorem~\ref{ssi} is given by
\begin{equation} \label{dd1}
 d_1 = \bd_1^{(r_1)} \star \text{Id}_2 \pm \text{Id}_1 \star
\bd_2^{(r_2)}\pm
 d_{\b_1} \star \text{Id}_2 \pm \text{Id}_1 \star d_{\b_2}
+ \d_{\b_1} \star \text{Id}_2 \pm \text{Id}_1 \star \d_{\b_2},
\end{equation}
where $\b = \b_1 \Si^{n-1}(\b_2)$ and $\ov{\b} = K_1 \# K_2$.
\end{thm} \qed

\noindent{\bf Remark:}
The notation in Theorem~\ref{boundary} and the
determination of signs can best be
explained by two
simple examples:
\[ d_1 (\r_1 \star \r_2)_1 =
(\bd_1^{(r_1)} \r_1 \star \r_2)_1
+ (-1)^{\mu({\r_1})}(\r_1 \star \bd_2^{(r_2)} \r_2)_1, \]
\[ d_1 (\r_1 \star s_{K_2}) = \bd_1^{(r_1)}\r_1 \star s_{K_2} 
+ (-1)^{\mu({\r_1})} (\r_1 \star \d_{\b_2} s_{K_2})_0 . \]
Note that $\mu({\r_1}^{(r_1)}) \equiv \mu({\r_1}) \pmod {2N}$.
We have extended our notation here in the obvious way: $(\sum m_j
\r_{1,j})\star \r_2 = \sum m_j (\r_{1, j} \star \r_2)$.
Theorem~\ref{boundary}
gives a full description of the differential $d_1$ in terms of two
special boundary maps and the $\Z$-graded boundary maps of
$\b_1$ and $\b_2$.
The reducible representations $s_{K_1}$ and $s_{K_2}$ do contribute
via the special maps $d_{\b_j}$ and $\d_{\b_j}$ ($j = 1, 2$)
in Definition~\ref{sbd}.
 
The differential $d_1: E^1_{p,q} \to E^1_{p-1, q}$ is possibly nonzero for
$q=0$ and $q=1$. Note that $E^1_{p,q}$ is generated from
$C_p^{(r_1, r_2)}(\p_{\b})$ by Lemma~\ref{generate}. Hence the differential
$d_1$ is calculated by the following: 
\begin{equation} \label{ri}
d_1(\r_1\star \r_2)_i = \sum {\hM}_J^1((\r_1\star \r_2)_i, (\r_1^{'}\star \r_2^{'})_i)
\cdot (\r_1^{'}\star \r_2^{'})_i, \ \ \ \text{for $i=0, 1$}; \end{equation}
\begin{equation} \label{rs}
d_1(\r_1\star s_{K_2}) = \sum {\hM}_J^1
(\r_1\star s_{K_2}, (\r_1^{'}\star \r_2^{'})_0) \cdot
(\r_1^{'}\star \r_2^{'})_0,
\end{equation}
for the critical points $(\r_1^{'} \star \r_2^{'})_0$ and
$(\r_1^{'} \star \r_2^{'})_1$ of the $U(1) = \r_1^{'} \star \r_2^{'}$
with Morse index $0$ and $1$. 
By (\ref{ri}) and (\ref{rs}), we have
\begin{equation} \label{rimas}
\mu_K((\r_1\star \r_2)_i) - \mu_K((\r_1^{'}\star \r_2^{'})_i) = 1; 
\end{equation}
\begin{equation} \label{rsmas}
\mu_K(\r_1\star s_2) - \mu_K ((\r_1^{'}\star \r_2^{'})_0) = 1.
\end{equation}
By Proposition~\ref{masl}, Corollary~\ref{maslc} and Lemma~\ref{generate}, 
\[\mu_K((\r_1\star \r_2)_i)=\mu_{K_1}(\r_1^{(r_1)}) + \mu_{K_2}
(\r_2^{(r_2)}) + i, \ \ \ i = 0, 1.\]
So the equation (\ref{rimas}) is reduced to
\[\mu_K((\r_1\star \r_2), (\r_1^{'}\star \r_2)) + \mu_K((\r_1^{'}\star \r_2),
(\r_1^{'}\star \r_2^{'})) = 1,\]
by Proposition~\ref{masl} and the path additivity of the Maslov 
index in \cite{clm1}. By the genericity, the moduli space of 
$J$-holomorphic curves is empty if the Maslov index is 
negative over the compact
monotone symplectic manifolds (see \cite{fl, lo, sz}). 
So we have the following two cases: 
\[\mu_K((\r_1\star \r_2), (\r_1^{'}\star \r_2)) =1, \ \ \ \
\mu_K((\r_1^{'}\star \r_2), (\r_1^{'}\star \r_2^{'})) = 0,\]
\[\mu_K((\r_1\star \r_2), (\r_1^{'}\star \r_2)) =0, \ \ \ \
\mu_K((\r_1^{'}\star \r_2), (\r_1^{'}\star \r_2^{'})) = 1.\]
For the case $\mu_K((\r_1^{'}\star \r_2), (\r_1^{'}\star \r_2^{'})) = 0$, we 
have $\mu_{K_2}(\r_2, \r_2^{'}) = 0$ over 
${\CR}^*(S^2 \setminus K_2)^{[i]}$, and
$\r_2 = \r_2^{'}$ by the monotonicity.
From the Maslov index calculation, the difference of Maslov index 1 must
have the difference of Malsov index 1 on one side and 0 on the other side.
Since both symplectic manifolds in the consideration are monotone, so we must
have a constant on the one side from the index calculation.

Now we study the moduli space of ${\CM}_J(\r_1\star \r_2, \r_1^{'}\star \r_2)$
for the Maslov index $\mu_{K_1}(\r_1, \r_1^{'}) = 1$. Since the boundary map
only counts the 1-dimensional moduli space, not $1 \pmod {2N(K)}$ dimensional
moduli space. For $u \in {\CM}^1_J(\r_1\star \r_2, \r_1^{'}\star \r_2)$,
we have $\mu_K(u) = 1$.
We are going to use a cobordism technique to relate the moduli space
${\CM}_J^1(\r_1\star \r_2, \r_1^{'}\star \r_2)$ with the
moduli space ${\CM}_J^1(\r_1, \r_1^{'}) \x \{\r_2\}$.

Let $J$ be an almost complex structure on ${\CR}^*(S^2 \setminus K)^{[i]}$
which is compatible with the symplectic structure. So $J$ is constructed
from an almost complex structure $J_j$ on
${\CR}^*(S^2 \setminus K_j)^{[i]}$ ($j=1, 2$). This gives the usual
decomposition into $J$-linear and $J$-anti-linear maps
\begin{equation} \label{jsplit}
\begin{split}
{\CE}_u & = \text{Hom}(T{\C}, u^*T{\CR}^*(S^2 \setminus K)^{[i]}) \\
& = T^{1,0}{\C}\ox_J u^*T{\CR}^*(S^2 \setminus K)^{[i]} \op
T^{0,1}{\C}\ox_J u^*T{\CR}^*(S^2 \setminus K)^{[i]}.\\
\end{split} \end{equation}
The subbundle ${\CE}^{0,1}$ is given by its fibre
\[{\CE}^{0,1}_u=L^p(T^{0,1}{\C}\ox_J u^*T{\CR}^*(S^2 \setminus K)^{[i]})
\to L^p_1(\text{Map}({\C}, {\CR}^*(S^2 \setminus K)^{[i]})).\]
Choose a trivialization of ${\CE}$ which is compatible with $J$. The operation
of $J$ on ${\CE}$ appears as an endomorphism of 
$L^p({\C}, \text{Hom}(T{\C}, u^*T{\CR}^*(S^2 \setminus K)^{[i]}))$.
The section $u \mapsto \pj u = \frac{1}{2}(du + J\circ du \circ i)$
is the projection of the section $u \mapsto du$ onto the subbundle
${\CE}^{0,1}$. The corresponding perturbed section 
$\ob_{J,H}$ is smooth, and its zero set is the moduli space
${\CM}_J(\r_1 \star \r_2, \r_1^{'} \star \r_2)$, where
$\ob_{J,H}u = \frac{\bd u}{\bd t} + J(\frac{\bd u}{\bd s} - X_s(u)) = 0$ and
$X_s$ is the corresponding Hamiltonian 
vector field of the Hamiltonian function $H:
{\CR}^*(S^2 \setminus K)^{[i]} \x {\R} \to {\R}$.

Note that the monotone symplectic structure on ${\CR}^*(S^2 \setminus K)^{[i]}$
is induced from the one on $Q_{n+m} = Q_n \x Q_m$ since $\b_1 \in B_n$
and $\b_2 \in B_m$ (see Lemma 2.2 in \cite{li}). The symplectic structure
$\o$ restricts on the symplectic structure $\o_j$ on
${\CR}^*(S^2 \setminus K_j)^{[i]}$ ($j=1, 2$). For the Hamiltonian function 
$H: {\CR}^*(S^2 \setminus K)^{[i]} \x {\R} \to {\R}$, we have
$\o (X_H, \cdot ) = dH$ and $\o_j (X_H, \cdot) = 
dH|_{{\CR}^*(S^2 \setminus K_j)^{[i]}}$.
Thus the restriction $H_j$ of $H$ on ${\CR}^*(S^2 \setminus K_j)^{[i]}$
is also a Hamiltonian function with respect to the restricted
symplectic structure $\o_j$.
We have $H_1(\p_{\b}(\mu_1), \cdot) = H_2(\p_{\b}(\mu_2), \cdot)$
and $H_j(x, s) = H_j(\p_{\b_j}(x), s+1)$. Both $H_1$ and $H_2$ agree
smoothly on the
identified meridian $\mu_1$ and $\mu_2$. So we denote by
$H = H_1 \star H_2$. Every Hamiltonian function
on ${\CR}^*(S^2 \setminus K)^{[i]}$ can be expressed in such a way.

Recall that the group operation on the space of Hamiltonians is given by
$H \circ H^{'} = H(x,s) + H^{'}((\p_H^s)^{-1}(x), s)$ and the inverse
$H^{-1}$ of $H$ under this operation is
$H^{-1}(x, s) = - H((\p_H^s)(x), s)$. So $H \circ H^{-1} = 0$.
Such an operation corresponds to the composition of the induced symplectic
diffeomorphisms on the symplectic manifold. We need to construct a
cobordism between the moduli spaces with fixed asymptotic values. This is
why we use this operation $H \circ H^{'}$ in our construction.

\begin{lm} \label{muf}
Over a monotone symplectic manifold $(M, \o)$, the zero dimensional
moduli space ${\CM}^0_J(x, x)$ of $J$-holomorphic curves consists of an
isolated constant flow $\{u(t) = x\}$.
\end{lm}
Proof: By the monotonicity, the symplectic action of any $J$-holomorphic
curves $u \in {\CM}^0_J(x, x)$ is zero. There are no nonconstant
$J$-holomorphic curves in $(M, \o)$. \qed

For the zero dimensional moduli space ${\CM}^0_{J_2, H_2}(\r_2, \r_2)$ on the 
monotone symplectic manifold $({\CR}^*(S^2 \setminus K_2)^{[i]}, \o_2)$,
we have an one parameter family of $J$-holomorphic curves
\[\ob_{J_2, H_{2, \t}} u = \frac{\bd u}{\bd t} + J_2(\frac{\bd u}{\bd s}
- X_{2, s}^{\t}(u)) = 0,\]
where $H_{2, \t} = H_2 \circ \t H^{-1}_2$ for $\t \in [0,1]$, and
$X_{2, s}^{\t}$ is the corresponding Hamiltonian vector field.
Such a family preserves the fixed point $\r_2$ at end, and at $\t = 0$ gives
the element in ${\CM}^0_{J_2, H_2}(\r_2, \r_2)$, and
at $\t = 1$ produces a unique element in ${\CM}^0_{J_2}(\r_2, \r_2) = 
\{ \r_2\}$ by Lemma~\ref{muf}. So the moduli space
$\ob_{J_2, H_{2, \t}}^{-1}(0)$ is regular at ends $\t = 0$ and
$\t = 1$, and the projection 
$\pi: \cup_{\t \in [0, 1]}\ob_{J_2, H_{2, \t}}^{-1}(0) \to [0, 1]$ 
is transverse at $\t$. So 
\[\# {\CM}^0_{J_2, H_2}(\r_2, \r_2) = \# {\CM}^0_{J_2, 0}(\r_2, \r_2)
= \pm 1, \]
by the unique constant trajectory flow $\{\r_2 \}$ in Lemma~\ref{muf}.

\begin{pro} \label{10c}
For $\r_1 \star \r_2$ and $\r_1^{'} \star \r_2$ in 
${\CR}^*(S^3 \setminus K)^{[i]}$ with $\mu_K(\r_1 \star \r_2,
\r_1^{'} \star \r_2) \equiv 1 \pmod {2N(K)}$, the one dimensional
moduli space ${\CM}^1_{J, H}(\r_1 \star \r_2, \r_1^{'} \star \r_2)$ has the
following property:
\[\# {\hM}^1_{J, H}(\r_1 \star \r_2, \r_1^{'} \star \r_2) =
\# {\hM}^1_{J_1, H_1}(\r_1, \r_1^{'}),\]
where $\# {\hM}^1_{J, H}(\r_1 \star \r_2, \r_1^{'} \star \r_2)$ is the algebraic
number of the zero dimensional manifold
${\hM}^1_{J, H}(\r_1 \star \r_2, \r_1^{'} \star \r_2)$.
\end{pro}
Proof: We consider the bundle of Banach spaces
${\CE}^{0,1} \to L^p_1(\text{Map}({\C}, {\CR}^*(S^2 \setminus K)^{[i]})) = B$
and the section $\ob_{J, H}$. We take the balanced maps for elements in $B$ with
asymptotic
values $\r_1 \star \r_2$ and $\r_1^{'} \star \r_2$.
Denote by $\hat{B}(\r_1 \star \r_2, \r_1^{'} \star \r_2)$.
The section $\ob_{J, H}$ is $t$-invariant and
$\ob_{J, H}$ induces a section on the quotient bundle
$\hat{\CE}^{0,1} \to \hat{B}(\r_1 \star \r_2, \r_1^{'} \star \r_2)$, where
$H = H_1 \star H_2$.
The fiber of $\hat{\CE}^{0,1}$ is given by
$L^p(T^{0,1}{\C} \ox_J u^*T{\CR}^*(S^2 \setminus K)^{[i]})$ with
$\int_{- \infty}^0|\n u|^2 = \int_0^{\infty }|\n u|^2$ for
$u \in \hat{B}(\r_1 \star \r_2, \r_1^{'} \star \r_2)$.
Again the zero set of $\ob_{J, H}$ is a zero dimensional manifold
${\hM}^1_{J, H}(\r_1 \star \r_2, \r_1^{'} \star \r_2)$ inside 
$\hat{B}(\r_1 \star \r_2, \r_1^{'} \star \r_2)$ for the one
dimensional component ${\hM}_{J, H}^1(\r_1 \star \r_2, 
\r_1^{'}\star \r_2)$. The orientation at 
each point in ${\hM}_{J, H}^1(\r_1 \star \r_2, \r_1^{'}\star \r_2)$
is the orientation on the trajectory flow used in defining the
symplectic Floer boundary. By the monotonicity and the Floer-Gromov
compactness,
\[\# {\hM}^1_{J, H}(\r_1 \star \r_2, \r_1^{'} \star \r_2) = \#
{\ob_{J, H}}^{-1}(0) \]
is well-defined.

Let $\t \in [0,1]$ and $\ob_{J, H_{\t}}$ be a section on
$\hat{\CE}^{0,1} \to \hat{B}(\r_1 \star \r_2, \r_1^{'} \star \r_2)$,
where $H_{\t} = H_1 \star (H_2 \circ \t H_2^{-1})$.
For each $\t$ the section $\ob_{J, H_{\t}}$ has the 
transverse zeros by \cite{fl}. For $\t = 0$, 
\[ \ob_{J, H_0}^{-1}(0) = \ob_{J,H_1 \star H_2}^{-1}(0) =
{\hM}^1_{J, H}(\r_1 \star \r_2, \r_1^{'} \star \r_2).\]
Let ${\hM}^1_{J, H_{\t}}(\r_1 \star \r_2, \r_1^{'} \star \r_2)$ be the zeros
of the section $\ob_{J, H_{\t}}$ for $\t \in [0, 1]$.
For $\t = 1$, the zero set $\ob_{J, H_1 \star (H_2 \circ H_2^{-1})}^{-1}(0)
= \ob_{J, H_1 \star 0}^{-1}(0)$ becomes solutions of
\begin{equation} \label{r1s}
\frac{\bd u}{\bd t} + J(\frac{\bd u}{\bd s} - X_{1, s}(u)) = 0.
\end{equation}
The $J$-holomorphic curve $u$ of (\ref{r1s}) satisfies the property that 
$u_1 \star u_2$ with $\mu_{K_1}(u_1) = 1$ and
$\mu_{K_2}(u_2) = 0$.
Thus $u_2 \in {\CM}^0_{J_2, H_2\circ H_2^{-1}}(\r_2, \r_2)$ is the
unique element $\{\r_2\}$ by Lemma~\ref{muf}.
Hence $u = u_1 \star \{\r_2\}$ can be identified as an element $u_1$ in
${\hM}^1_{J_1, H_1}(\r_1, \r_1^{'})$.

Then $\t = 0$ and $\t =1$ in $[0, 1]$ are regular values of the projection 
$\pi_1: \cup_{\t \in [0, 1]}{\hM}_{J,H_{\t}} \to [0, 1]$ since the 
one dimensional moduli space $\cup_{\t \in [0, 1]}{\hM}_{J,H_{\t}}$ is 
transverse to $\hat{B}(\r_1 \star \r_2, \r_1^{'} \star \r_2)$
(see \cite{dk, fl, lo, sz}). So the parameterized moduli space
$\cup_{\t \in [0, 1]}{\hM}_{J,H_{\t}}$ is an one dimensional submanifold
in the product space $[0,1] \x \hat{B}(\r_1 \star \r_2, \r_1^{'} \star \r_2)$
with oriented boundary components
$- {\hM}^1_{J, H}(\r_1 \star \r_2, \r_1^{'} \star \r_2)$
and ${\hM}^1_{J_1, H_1}(\r_1, \r_1^{'}) \star \{\r_2\}$.
Each boundary component is compact zero dimensional manifold and
$\cup_{\t \in [0, 1]}{\hM}_{J,H_{\t}}$ is also compact by the 
Floer-Gromov compactness theorem.
So ${\hM}^1_{J, H}(\r_1 \star \r_2, \r_1^{'} \star \r_2)$ and 
${\hM}^1_{J_1, H_1}(\r_1, \r_1^{'}) \star \{\r_2\}$ are oriented cobordant:
\begin{equation*} \begin{split}
0 & = \bd (\cup_{\t \in [0, 1]}{\hM}_{J,H_{\t}}) \\
 & = - \# {\ob_{J, H_{\t}=H_0}}^{-1}(0) + \# {\ob_{J, H_{\t} =H_0}}^{-1}(0) \\
 & = - \# {\hM}^1_{J, H}(\r_1 \star \r_2, \r_1^{'} \star \r_2)
+ \# {\hM}^1_{J_1, H_1}(\r_1, \r_1^{'}) \star \{\r_2\} \\
& = - \# {\hM}^1_{J, H}(\r_1 \star \r_2, \r_1^{'} \star \r_2)
+ \# {\hM}^1_{J_1, H_1}(\r_1, \r_1^{'}),
\end{split} \end{equation*}
since the orientation of $u_1 \star \{\r_2\}$ is compatible with the
orientations of $\r_1 \star \r_2$ and $\r_1^{'} \star \r_2$.
Hence the result follows. \qed

\noindent{\bf Proof of Theorem~\ref{boundary}:}
By Corollary~\ref{maslc} and Proposition~\ref{10c}, we have that 
all possible one dimensional moduli spaces
${\CM}^1_{J, H}(\r_1 \star \r_2, \r_1^{'} \star \r_2^{'})$ are given
by the asymptotic values same on one of
${\CR}^*(S^3 \setminus K_j)^{[i]} (j=1, 2)$
and differ by one Maslov index on the other. There is no nontrivial
$J$-holomorphic curve with negative Maslov index on
the monotone symplectic manifolds. 
So the first differential $d_1$ is given by
\[d_1 = \partial_1^{(r_1)} \star \text{Id}_2 \pm \text{Id}_1 \star
\bd_2^{(r_2)}\pm
 d_{\b_1} \star \text{Id}_2 \pm \text{Id}_1 \star d_{\b_2}
+ \d_{\b_1} \star \text{Id}_2 \pm \text{Id}_1 \star \d_{\b_2},\]
for all possible $\r_1 \star \r_2$ in Proposition~\ref{ge}.
One needs to check that (\ref{dd1}) is indeed a differential, i.e.,
$d_1 \circ d_1 = 0$. This follows by straightforward calculation by
using the remark after Theorem~\ref{boundary}. \qed

\noindent{\bf Remarks}: (i) The identification of $d_1$ is highly depending
on the monotonicity of ${\CR}^*(S^2 \setminus K)^{[i]}$.
By \cite{fl} and Proposition~\ref{masl}, we reduce the problem to the
moduli space with special asymptotic values (one of
$\r_j$ and $\r_j^{'}$ ($j=1, 2$) agrees). Then we construct a cobordism
to compute the first differential $d_1$.

(ii) In \cite{li1, li3}, we studied the structure of all moduli spaces
with index $\leq 4$ by a gluing result. Using certain cohomology classes, 
we identified all the higher differentials in \cite{li3}.
In order to identify $d_2$ in Theorem~\ref{ssi}, we need to get the local
structure of $J$-holomorphic curves on ${\CR}^*(S^2 \setminus K)^{[i]}$
which is similar to the gluing result in \cite{li1}, and then find the 
correct cohomology class of degree one for the differential $d_2$.
We will discuss this elsewhere. 

(iii) Note that ${\CR}^*(S^2 \setminus K)$ can be identified as a symplectic
fiber product of ${\CR}^*(S^2 \setminus K_1)$ and ${\CR}^*(S^2 \setminus K_2)$.
We hope that our discussion in this paper may shed a light on the study
of the symplectic Floer homology of the fiber product of two symplectic
manifolds.

\end{document}